\documentclass[12pt,reqno]{amsart}
\usepackage{amssymb}
\title[Families of hypersurfaces]
{Complex multiplication, Griffiths-Yukawa couplings, and rigidity for families of hypersurfaces}
\author[Eckart Viehweg]{Eckart Viehweg}
\address{Universit\"at Essen, FB6 Mathematik, 45117 Essen, Germany}
\email{viehweg@uni-essen.de}
\thanks{This work has been supported by the Institute of Mathematical Science at the Chinese
University of Hong Kong and by the ``DFG-Schwerpunktprogramm
Globale Methoden in der Komplexen Geometrie''. The second named author
is supported by grants from the Research Grants Council of the Hong Kong
Special Administrative Region, China (Project No. CUHK 4034/02P)
and from the Institute of  Mathematical Sciences at the Chinese University
of Hong Kong (Program in Algebraic Geometry)}
\author[Kang Zuo]{Kang Zuo}
\address{The Chinese University of Hong Kong, Department of Mathematics,
Shatin, Hong Kong}
\email{kzuo@math.cuhk.edu.hk}
\begin{document}
%%%%%%%%%%%%%%%%%%%% Text italic %%%%%%%%%%%%%%%%%%%%%%%%%%%%
\theoremstyle{plain}
\newtheorem{thm}{Theorem}[section]
\newtheorem{theorem}[thm]{Theorem}
\newtheorem{lemma}[thm]{Lemma}
\newtheorem{corollary}[thm]{Corollary}
\newtheorem{proposition}[thm]{Proposition}
\newtheorem{construction}[thm]{Construction}
%%%%%%%%%%%%%%%%%%%% Text roman %%%%%%%%%%%%%%%%%%%%%%%%%%%%%
\theoremstyle{definition}
\newtheorem{notations}[thm]{Notations}
\newtheorem{notation}[thm]{Notation}
\newtheorem{problem}[thm]{Problem}
\newtheorem{remark}[thm]{Remark}
\newtheorem{remarks}[thm]{Remarks}
\newtheorem{definition}[thm]{Definition}
\newtheorem{claim}[thm]{Claim}
\newtheorem{assumption}[thm]{Assumption}
\newtheorem{assumptions}[thm]{Assumptions}
\newtheorem{properties}[thm]{Properties}
\newtheorem{example}[thm]{Example}
\newtheorem{conjecture}[thm]{Conjecture}
\numberwithin{equation}{section}
%%%%%%%% Diagram macros, etc. %%%%%%%%%%%%%%%%%%%%%%%%%%%%%%%
\catcode`\@=11
% General macros
\def\opn#1#2{\def#1{\mathop{\kern0pt\fam0#2}\nolimits}}
\def\bold#1{{\bf #1}}%
\def\underrightarrow{\mathpalette\underrightarrow@}
\def\underrightarrow@#1#2{\vtop{\ialign{$##$\cr
 \hfil#1#2\hfil\cr\noalign{\nointerlineskip}%
 #1{-}\mkern-6mu\cleaders\hbox{$#1\mkern-2mu{-}\mkern-2mu$}\hfill
 \mkern-6mu{\to}\cr}}}
\let\underarrow\underrightarrow
\def\underleftarrow{\mathpalette\underleftarrow@}
\def\underleftarrow@#1#2{\vtop{\ialign{$##$\cr
 \hfil#1#2\hfil\cr\noalign{\nointerlineskip}#1{\leftarrow}\mkern-6mu
 \cleaders\hbox{$#1\mkern-2mu{-}\mkern-2mu$}\hfill
 \mkern-6mu{-}\cr}}}
% Rectangular Commutative diagrams
\let\amp@rs@nd@\relax
\newdimen\ex@
\ex@.2326ex
\newdimen\bigaw@
\newdimen\minaw@
\minaw@16.08739\ex@
\newdimen\minCDaw@
\minCDaw@2.5pc
\newif\ifCD@
\def\minCDarrowwidth#1{\minCDaw@#1}
\newenvironment{CD}{\@CD}{\@endCD}
\def\@CD{\def\A##1A##2A{\llap{$\vcenter{\hbox
 {$\scriptstyle##1$}}$}\Big\uparrow\rlap{$\vcenter{\hbox{%
$\scriptstyle##2$}}$}&&}%
\def\V##1V##2V{\llap{$\vcenter{\hbox
 {$\scriptstyle##1$}}$}\Big\downarrow\rlap{$\vcenter{\hbox{%
$\scriptstyle##2$}}$}&&}%
\def\={&\hskip.5em\mathrel
 {\vbox{\hrule width\minCDaw@\vskip3\ex@\hrule width
 \minCDaw@}}\hskip.5em&}%
\def\verteq{\Big\Vert&&}%
\def\noarr{&&}%
\def\vspace##1{\noalign{\vskip##1\relax}}\relax\iffalse{%
\fi\let\amp@rs@nd@&\iffalse}\fi
 \CD@true\vcenter\bgroup\relax\iffalse{%
\fi\let\\=\cr\iffalse}\fi\tabskip\z@skip\baselineskip20\ex@
 \lineskip3\ex@\lineskiplimit3\ex@\halign\bgroup
 &\hfill$\m@th##$\hfill\cr}
\def\@endCD{\cr\egroup\egroup}
% Horizontal arrows with "sliding" length
\def\>#1>#2>{\amp@rs@nd@\setbox\z@\hbox{$\scriptstyle
 \;{#1}\;\;$}\setbox\@ne\hbox{$\scriptstyle\;{#2}\;\;$}\setbox\tw@
 \hbox{$#2$}\ifCD@
 \global\bigaw@\minCDaw@\else\global\bigaw@\minaw@\fi
 \ifdim\wd\z@>\bigaw@\global\bigaw@\wd\z@\fi
 \ifdim\wd\@ne>\bigaw@\global\bigaw@\wd\@ne\fi
 \ifCD@\hskip.5em\fi
 \ifdim\wd\tw@>\z@
 \mathrel{\mathop{\hbox to\bigaw@{\rightarrowfill}}\limits^{#1}_{#2}}\else
 \mathrel{\mathop{\hbox to\bigaw@{\rightarrowfill}}\limits^{#1}}\fi
 \ifCD@\hskip.5em\fi\amp@rs@nd@}
\def\<#1<#2<{\amp@rs@nd@\setbox\z@\hbox{$\scriptstyle
 \;\;{#1}\;$}\setbox\@ne\hbox{$\scriptstyle\;\;{#2}\;$}\setbox\tw@
 \hbox{$#2$}\ifCD@
 \global\bigaw@\minCDaw@\else\global\bigaw@\minaw@\fi
 \ifdim\wd\z@>\bigaw@\global\bigaw@\wd\z@\fi
 \ifdim\wd\@ne>\bigaw@\global\bigaw@\wd\@ne\fi
 \ifCD@\hskip.5em\fi
 \ifdim\wd\tw@>\z@
 \mathrel{\mathop{\hbox to\bigaw@{\leftarrowfill}}\limits^{#1}_{#2}}\else
 \mathrel{\mathop{\hbox to\bigaw@{\leftarrowfill}}\limits^{#1}}\fi
 \ifCD@\hskip.5em\fi\amp@rs@nd@}
% Rectangular commutative diagrams with diagonal arrows
\newenvironment{CDS}{\@CDS}{\@endCDS}
\def\@CDS{\def\A##1A##2A{\llap{$\vcenter{\hbox
 {$\scriptstyle##1$}}$}\Big\uparrow\rlap{$\vcenter{\hbox{%
$\scriptstyle##2$}}$}&}%
\def\V##1V##2V{\llap{$\vcenter{\hbox
 {$\scriptstyle##1$}}$}\Big\downarrow\rlap{$\vcenter{\hbox{%
$\scriptstyle##2$}}$}&}%
\def\={&\hskip.5em\mathrel
 {\vbox{\hrule width\minCDaw@\vskip3\ex@\hrule width
 \minCDaw@}}\hskip.5em&}
\def\verteq{\Big\Vert&}
\def\novarr{&}
\def\noharr{&&}
\def\SE##1E##2E{\slantedarrow(0,18)(4,-3){##1}{##2}&}
\def\SW##1W##2W{\slantedarrow(24,18)(-4,-3){##1}{##2}&}
\def\NE##1E##2E{\slantedarrow(0,0)(4,3){##1}{##2}&}
\def\NW##1W##2W{\slantedarrow(24,0)(-4,3){##1}{##2}&}
\def\slantedarrow(##1)(##2)##3##4{%
\thinlines\unitlength1pt\lower 6.5pt\hbox{\begin{picture}(24,18)%
\put(##1){\vector(##2){24}}%
\put(0,8){$\scriptstyle##3$}%
\put(20,8){$\scriptstyle##4$}%
\end{picture}}}
\def\vspace##1{\noalign{\vskip##1\relax}}\relax\iffalse{%
\fi\let\amp@rs@nd@&\iffalse}\fi
 \CD@true\vcenter\bgroup\relax\iffalse{%
\fi\let\\=\cr\iffalse}\fi\tabskip\z@skip\baselineskip20\ex@
 \lineskip3\ex@\lineskiplimit3\ex@\halign\bgroup
 &\hfill$\m@th##$\hfill\cr}
\def\@endCDS{\cr\egroup\egroup}
% Triangular commutative diagrams
\newdimen\TriCDarrw@
\newif\ifTriV@
\newenvironment{TriCDV}{\@TriCDV}{\@endTriCD}
\newenvironment{TriCDA}{\@TriCDA}{\@endTriCD}
\def\@TriCDV{\TriV@true\def\TriCDpos@{6}\@TriCD}
\def\@TriCDA{\TriV@false\def\TriCDpos@{10}\@TriCD}
\def\@TriCD#1#2#3#4#5#6{%
\setbox0\hbox{$\ifTriV@#6\else#1\fi$}
\TriCDarrw@=\wd0 \advance\TriCDarrw@ 24pt
\advance\TriCDarrw@ -1em
\def\SE##1E##2E{\slantedarrow(0,18)(2,-3){##1}{##2}&}
\def\SW##1W##2W{\slantedarrow(12,18)(-2,-3){##1}{##2}&}
\def\NE##1E##2E{\slantedarrow(0,0)(2,3){##1}{##2}&}
\def\NW##1W##2W{\slantedarrow(12,0)(-2,3){##1}{##2}&}
\def\slantedarrow(##1)(##2)##3##4{\thinlines\unitlength1pt
\lower 6.5pt\hbox{\begin{picture}(12,18)%
\put(##1){\vector(##2){12}}%
\put(-4,\TriCDpos@){$\scriptstyle##3$}%
\put(12,\TriCDpos@){$\scriptstyle##4$}%
\end{picture}}}
\def\={\mathrel {\vbox{\hrule
   width\TriCDarrw@\vskip3\ex@\hrule width
   \TriCDarrw@}}}
\def\>##1>>{\setbox\z@\hbox{$\scriptstyle
 \;{##1}\;\;$}\global\bigaw@\TriCDarrw@
 \ifdim\wd\z@>\bigaw@\global\bigaw@\wd\z@\fi
 \hskip.5em
 \mathrel{\mathop{\hbox to \TriCDarrw@
{\rightarrowfill}}\limits^{##1}}
 \hskip.5em}
\def\<##1<<{\setbox\z@\hbox{$\scriptstyle
 \;{##1}\;\;$}\global\bigaw@\TriCDarrw@
 \ifdim\wd\z@>\bigaw@\global\bigaw@\wd\z@\fi
 \mathrel{\mathop{\hbox to\bigaw@{\leftarrowfill}}\limits^{##1}}
 }
 \CD@true\vcenter\bgroup\relax\iffalse{\fi\let\\=\cr\iffalse}\fi
 \tabskip\z@skip\baselineskip20\ex@
 \lineskip3\ex@\lineskiplimit3\ex@
 \ifTriV@
 \halign\bgroup
 &\hfill$\m@th##$\hfill\cr
#1&\multispan3\hfill$#2$\hfill&#3\\
&#4&#5\\
&&#6\cr\egroup%
\else
 \halign\bgroup
 &\hfill$\m@th##$\hfill\cr
&&#1\\%
&#2&#3\\
#4&\multispan3\hfill$#5$\hfill&#6\cr\egroup
\fi}
\def\@endTriCD{\egroup}
%%%%%%%%%%%%%%%  End of diagram macros.  %%%%%%%%%%%%%%%%%%%%%%%%%
% Skriptbuchstaben
\newcommand{\sA}{{\mathcal A}}
\newcommand{\sB}{{\mathcal B}}
\newcommand{\sC}{{\mathcal C}}
\newcommand{\sD}{{\mathcal D}}
\newcommand{\sE}{{\mathcal E}}
\newcommand{\sF}{{\mathcal F}}
\newcommand{\sG}{{\mathcal G}}
\newcommand{\sH}{{\mathcal H}}
\newcommand{\sI}{{\mathcal I}}
\newcommand{\sJ}{{\mathcal J}}
\newcommand{\sK}{{\mathcal K}}
\newcommand{\sL}{{\mathcal L}}
\newcommand{\sM}{{\mathcal M}}
\newcommand{\sN}{{\mathcal N}}
\newcommand{\sO}{{\mathcal O}}
\newcommand{\sP}{{\mathcal P}}
\newcommand{\sQ}{{\mathcal Q}}
\newcommand{\sR}{{\mathcal R}}
\newcommand{\sS}{{\mathcal S}}
\newcommand{\sT}{{\mathcal T}}
\newcommand{\sU}{{\mathcal U}}
\newcommand{\sV}{{\mathcal V}}
\newcommand{\sW}{{\mathcal W}}
\newcommand{\sX}{{\mathcal X}}
\newcommand{\sY}{{\mathcal Y}}
\newcommand{\sZ}{{\mathcal Z}}
% Sonderbuchstaben mit Doppellinie
\newcommand{\A}{{\mathbb A}}
\newcommand{\B}{{\mathbb B}}
\newcommand{\C}{{\mathbb C}}
\newcommand{\D}{{\mathbb D}}
\newcommand{\E}{{\mathbb E}}
\newcommand{\F}{{\mathbb F}}
\newcommand{\G}{{\mathbb G}}
\newcommand{\BH}{{\mathbb H}}
\newcommand{\I}{{\mathbb I}}
\newcommand{\J}{{\mathbb J}}
\newcommand{\LL}{{\mathbb L}}
\newcommand{\M}{{\mathbb M}}
\newcommand{\N}{{\mathbb N}}
\newcommand{\BP}{{\mathbb P}}
\newcommand{\Q}{{\mathbb Q}}
\newcommand{\R}{{\mathbb R}}
\newcommand{\BS}{{\mathbb S}}
\newcommand{\T}{{\mathbb T}}
\newcommand{\U}{{\mathbb U}}
\newcommand{\V}{{\mathbb V}}
\newcommand{\W}{{\mathbb W}}
\newcommand{\X}{{\mathbb X}}
\newcommand{\Y}{{\mathbb Y}}
\newcommand{\Z}{{\mathbb Z}}
\newcommand{\id}{{\rm id}}
\newcommand{\rk}{{\rm rank}}
\newcommand{\roundup}[1]{\ulcorner{#1}\urcorner}
\newcommand{\rounddown}[1]{\llcorner{#1}\lrcorner}
\newcommand{\Hg}{{\rm Hg}}
\newcommand{\Sl}{{\rm Sl}}
\newcommand{\Gl}{{\rm Gl}}
%%%%%%%%%%%%%%%%%%%%%%%%%%%%%%%%%%%%%%%%%%%%%%%%%%%%%%%%%%%%%%
\begin{abstract}
Let $\sM_{d,n}$ be the moduli stack of hypersurfaces $X\subset \BP^n$ of degree $d\geq n+1$,
and let $\sM_{d,n}^{(1)}$ be the sub-stack, parameterizing hypersurfaces obtained
as a $d$ fold cyclic covering of $\BP^{n-1}$ ramified over a hypersurface of degree $d$.
Iterating this construction, one obtains $\sM_{d,n}^{(\nu)}$.

We show that $\sM_{d,n}^{(1)}$ is rigid in $\sM_{d,n}$, although for $d<2n$ the Griffiths-Yukawa
coupling degenerates. However, for all $d\geq n+1$ the sub-stack $\sM^{(2)}_{d,n}$ deforms.

We calculate the exact length of the Griffiths-Yukawa coupling
over $\sM_{d,n}^{(\nu)}$, and we construct a $4$-dimensional
family of quintic hypersurfaces $g:\sZ\to T$ in $\BP^4$, and a
dense set of points $t$ in $T$, such that $g^{-1}(t)$ has complex
multiplication.
\end{abstract}
\maketitle
%%%%%%%%%%% Introduction %%%%%%%%%%%%%%%%%%%%%%%%%%%%%%%%%%%
\section*{Introduction}
Let $\sM_{d,n}$ denote the moduli stack of hypersurfaces of degree $d\geq 2$
in the complex projective space $\BP^{n}$, and let $M_{d,n}$ be
the corresponding coarse moduli scheme.
Hence $\sM_{d,n}(\C)$ classifies for $n>1$ pairs $(X,\sL)$ with $X$ a nonsingular
manifold of dimension $n-1$, and with $\sL$ a very ample invertible sheaf
with $h^0(X,\sL)=n+1$. We will frequently write $\sO_X(1)$ instead of $\sL$.

A morphism $S\to M_{d,n}$ factors through the moduli stack
$\sM_{d,n}$ if it is induced by some pair $(f: \sX \to S, \sL)\in \sM_{d,n}(S)$.
Then $R^if_*\sL$ is zero for $i>0$, and locally free of rank $n+1$, for
$n>1$ and $i=0$. Moreover $\sX$ is embedded over $S$ in $\BP(f_*\sL)$.
We call $(f: \sX \to S, \sL)\in \sM_{d,n}(S)$ a universal family, whenever
the induced morphism $S \to M_{d,n}$ is dominant and generically finite
(and if $(f: \sX \to S, \sL)$ is normalized, as defined in \ref{1.2}).

It is the aim of this article, to study certain sub-moduli stacks $\sM^{(\nu)}_{d,n}$
of $\sM_{d,n}$, for $d\geq n+1$. Roughly speaking, starting from a family
$$(f: \sX \to S, \sL)\in \sM_{d,n-1}(S)$$
one can construct a new family in $\sM_{d,n}(S)$ by taking a cyclic covering of degree $d$
of the projective bundle $\BP(f_*\sL)$, totally ramified over the divisor $\sX$ and nowhere else
(see Section \ref{cyclic}).
The moduli stack of such families will be denoted by $\sM^{(1)}_{d,n}$.
Repeating this process $\nu$ times, starting of course with families in $\sM_{d,n-\nu}(S)$
one obtains the families in $\sM^{(\nu)}_{d,n}(S)$. The corresponding coarse moduli scheme
will be denoted by $M^{(\nu)}_{d,n}$.

Given $(f: \sX \to S, \sL)\in \sM_{d,n}(S)$ consider the variation of polarized
Hodge structures $R^{n-1}f_*\Q_{\sX}$, or the corresponding Higgs bundle
(called system of Hodge bundles by Simpson)
$$
\big(E=\bigoplus_{p+q=n-1} E^{p,q}, \ \theta=\bigoplus_{p+q=n-1}\theta_{p,q}\big),
$$
where $E^{p,q}=R^qf_*\Omega^p_{\sX/S}$ and where the Higgs field
$$\theta_{p,q}:E^{p,q} \to E^{p-1,q+1}\otimes \Omega_S^1$$
is given by the cup product with the Kodaira Spencer map (see \cite{VZ3}, for example), i.e. by the
edge morphisms of the wedge products of the tautological exact sequence
$$
0\>>> f^*\Omega^1_S \>>> \Omega^1_\sX \>>> \Omega^1_{\sX/S} \>>> 0.
$$
The $i$-th iterated cup product with the Kodaira Spencer map defines a morphism
\begin{multline*}
\theta^i:E^{0,n-1} \>\theta_{0,n-1}>> E^{1,n-2}\otimes \Omega^1_S
\>\theta_{1,n-2}>>\\ E^{2,n-3}\otimes S^2(\Omega^1_S) \>\theta_{2,n-3}>> \cdots
\>\theta_{i,n-i}>> E^{i,n-1-i}\otimes S^i(\Omega^1_S).
\end{multline*}
For $i=n-1$ one obtains a coupling
$$
\theta^{n-1}:E^{0,n-1} \>>> E^{n-1,0}\otimes S^{n-1}(\Omega^1_S)={E^{0,n-1}}^\vee
\otimes S^{n-1}(\Omega^1_S),
$$
which for families of hypersurfaces has been studied by Carlson, Green and Griffiths
(see \cite{CG}). For families of Calabi-Yau manifolds, the importance of this
coupling was brought up by physicists and they studied it in detail.
We will call it the Yukawa coupling, if the fibres are Calabi-Yau manifolds,
and the Griffiths-Yukawa coupling in general.
We define its length to be
$$
\varsigma(f)={\rm Min}\{i\geq 1 ; \ \theta^i=0\}-1.
$$
We will write $\varsigma(\sM^{(\nu)}_{d,n})$ instead of $\varsigma(f)$,
if the family is universal, hence if the induced morphism $S\to M^{(\nu)}_{d,n}$ is
dominant and generically finite.
\begin{theorem}[Section \ref{cyclic} and \ref{variation}]\label{6.5} For $n\geq 3$ and $d\geq n+1$
consider the sub-moduli stacks
$$
\sM^{(n-1)}_{d,n} \subset \sM_{d,n}^{(n-2)} \subset \cdots \subset \sM_{d,n}^{(2)}\subset
\sM_{d,n}^{(1)} \subset \sM_{d,n}.
$$
Then
\begin{enumerate}
\item[(a)] $\varsigma(\sM^{(\ell)}_{d,n})=n-\ell$ for $\ell = n-[\frac{d}{2}]+1, \ldots ,n-1$.
\item[(b)] $\varsigma(\sM^{(\ell)}_{d,n})=n-\ell-1$ for $\ell = 1,\ldots , n-[\frac{d}{2}]$.
\end{enumerate}
\end{theorem}

Remark that for $\ell=1$ Theorem \ref{6.5} implies that $\varsigma(\sM^{(1)}_{d,n})=n-2$,
if and only if $2n>d$. In particular, for families of Calabi-Yau hypersurfaces
belonging to $\sM^{(1)}_{n+1,n}$, the Yukawa coupling is always zero.

For families of canonically polarized manifolds, or for families of minimal models
of Kodaira dimension zero, the maximality of $\varsigma(f)$ implies rigidity,
i.e. that for $\dim(S)>0$ and $\dim(T)>0$ there can not exist a generically
finite morphism from $S\times T$ to the corresponding moduli scheme, which is induced
by a family. We will say shortly, that $S$ is rigid in the moduli stack or, if
$(f:\sX \to S,\sL)$ is a universal family for a sub-moduli stack $\sM$ of
$\sM_{d,n}$, that $\sM$ is rigid in $\sM_{d,n}$.

The observation that the maximality of the length of the Griffiths-Yukawa coupling implies
rigidity has implicitly been used in \cite{VZ3}, Proof of 6.4 and 6.5,
and it was stated explicitly in the survey \cite{VZ4}, Section 8.
A similar result has been shown by S. Kov\'acs and, for families of Calabi-Yau
manifolds, by  K-F. Liu, A. Todorov, S.-T. Yau and the second named author
in \cite{LTYZ}.

Together with Theorem \ref{6.5} it implies that $\sM^{(1)}_{d,n}$
is rigid in $\sM_{d,n}$ for $d\geq 2n$. As we will see in Section \ref{iterated}
the same holds true for $n+1 \leq d < 2n$, although $\varsigma(f)<n-1$.

\begin{theorem}\label{8.3} A universal family $g:\sZ\to S$ for $\sM^{(1)}_{d,n}$ is rigid.
\end{theorem}
As well known (see \ref{1.3}), for $n\geq 4$ or for $n=3$ and $d\geq 5$
all deformations of a hypersurface in $\BP^n$ are again hypersurfaces in $\BP^n$,
hence the rigidity in Theorem \ref{8.3} is independent from the polarization
chosen for $g:\sZ \to S$.

Let us remark, that the rigidity of families also follow from a strong positivity
property of the sheaf of logarithmic differential forms on compactifications of
$\sM_{d,n}$. To be more precise, let $Y$ be a projective manifold and $S\subset Y$
the complement of a normal crossing divisor $\Gamma$. Assume that for some $(d,n)$
and for all generically finite morphisms $S\to M_{d,n}$, factoring through the moduli
stack, the sheaf $\Omega^1_Y(\log \Gamma)$ is big (see \cite{VZ3}, Definition 1.1).
Then all families $f:\sX\to S$, as above, are rigid. By \cite{ACT} the moduli stack
$M_{3,3}$ has this property. However, as we will see, this no longer holds true for
$n\geq 3$ and $d\geq n+1$.

\begin{theorem}\label{5.5} Assume that for $n\geq 3$, for $d\geq n+1$,
and for some $\nu \leq n-1$ the morphism $S_\nu \to M^{(\nu)}_{d,n}$
is generically finite and induced by a family. Then there exists a
$(d-3)$-dimensional manifold $T$, and for $r=[\frac{\nu}{2}]$ a generically
finite morphism $S_\nu \times T^{\times r} \to M_{d,n}$
which is induced by a family.
\end{theorem}
In particular, the moduli stack $\sM_{d,n}^{(2)}$ always deforms in a non-trivial way in
$\sM_{d,n}$.

As pointed out by S.-T. Yau, the sub-moduli stack $\sM^{(3)}_{5,4}$ has been studied before 
by S. Ferrara and J.~Louis in \cite{FL}. There it is shown, that $M^{(3)}_{5,4}$ has 
a natural structure of a ball quotient (see Remark \ref{7.7}), and that 
$\varsigma(\sM^{(3)}_{5,4})\leq 2$.

As a byproduct of the calculation of variations of Hodge structures for families in
$\sM_{d,n}^{(2)}$ (Section \ref{variation}, see also \cite{vGI}, Section 3) we show that for a
universal family $g:\sZ \to S$ in $\sM^{(3)}_{5,4}$ the set of CM-points is dense in $S$, i.e.
the set of points $s\in S$ where the fibre $g^{-1}(s)$ has complex multiplication (see Section
\ref{quintic}). Together with Theorem \ref{5.5} this will imply that the Zariski closure
of the set of CM-points contains a $4$-dimensional subvariety $M$.
\begin{theorem}\label{0.2}
There exists infinitely many quintic threefolds with complex multiplication.
More precisely, there exists a finite and rigid map
$$M_{5,1}\times M_{5,1}\>>> M_{5,4}$$
with image $M$, such that the CM-points are dense in $M$.
\end{theorem}
The arguments used in the proof of \ref{0.2} do not extend to the case $n\geq 5$ and $d=n+1$.
They are related to the ones used by B. van Geemen and E. Izadi in \cite{vGI}.
In particular in \cite{vGI} quintic surfaces with complex multiplication are studied in Section 6.\\

This note grew out of discussions started when the first named author visited
the Institute of Mathematical Science and the Department of Mathematics at the Chinese
University of Hong Kong. His contributions to the present version
were written during two visits to the I.H.E.S., Bures sur Yvette.
He would like to thank the members of those three Institutes for their
hospitality. 

Shing-Tung Yau, drew our attention to the the work of S. Ferrara and J.~Louis
\cite{FL}, an article he was studying to understand similar questions.
Bert van Geemen found an error in the first version of Section \ref{quintic},
and he pointed out the relation between our Sections \ref{variation} and
\ref{quintic} the Sections 3 and 6 of \cite{vGI}.
We both would like to thank them, and H\'el\`ene Esnault for their interest and help. 

\section{Moduli of Hypersurfaces and the Jacobian ring}

Let us recall some well known vanishing theorems for sheaves of (logarithmic) differential forms
on $\BP^n$.
\begin{lemma}\label{1.1}
Let $F$ be a non-singular hypersurface in $\BP^n$ of degree $d$.
\begin{enumerate}
\item[(a)] $H^q(\BP^n,\Omega^p_{\BP^n}\otimes \sO_{\BP^n}(\nu))=0$ for
\begin{enumerate}
\item[i.] $q=0$, and $\nu \leq p$.
\item[ii.] $0 < q$ and $\nu > 0$.
\item[iii.] $q<n$ and $\nu < 0$.
\item[iv.] $q=n$, and $\nu \geq p-n $.
\item[v.] $p\neq q$, and $\nu=0$.
\end{enumerate}
\item[(b)] $H^q(\BP^n,\Omega^p_{\BP^n}(\log F)\otimes \sO_{\BP^n}(\nu))=0$ for
\begin{enumerate}
\item[i.] $p+q < n$, and $\nu <0$.
\item[ii.] $p+q > n$, and $\nu > -d$.
\end{enumerate}
\end{enumerate}
\end{lemma}
\begin{proof}
a) In \cite{Laz}, 7.3.9, for example, one finds a proof for i) and for ii).
Then Serre duality implies iii) and iv), and v) is obvious.

b, i) is a very special and known case of \cite{EV1}, 6.4., and ii) follows
again by Serre duality.
\end{proof}
Consider for $n>1$ and $d\geq n+1$ a family $(f: \sX \to S,
\sL)\in \sM_{d,n}(S)$ of hypersurfaces in $\BP^n$ of degree $d$.
The polarization $\sL$ on $\sX$ is only determined up to $\otimes
f^*\sB$ for invertible sheaves $\sB$ on $S$. For $d=n+1$ the sheaf
$\omega_\sX$ is trivial, and for $d>n+1$ it is equivalent to
$\sL^{d-n-1}$. So for $d>n+1$ the moduli stack $\sM_{d,n}$ is
finite over the moduli stack of canonically polarized manifolds.

\begin{notation}\label{1.2}
We will call $(f: \sX \to S, \sL)\in \sM_{d,n}(S)$ normalized, if
$$\sO_{\BP(f_*\sL)}(\sX) = \sO_{\BP(f_*\sL)}(d),$$
and we will write $\sE=f_*\sL$. If $\sM$ is a sub-moduli stack with coarse moduli scheme
$M$ we call $(f: \sX \to S, \sL)$ a universal family for $\sM$, if it is normalized, and
if the induced morphism $S\to M$ is dominant and generically finite.
Often we will not mention the polarization, and just write $(f: \sX \to S)\in \sM_{d,n}(S)$.
\end{notation}
Replacing $S$ by some finite covering, one can always choose $\sL$
such that the family is normalized.
\begin{lemma}\label{1.3}
Assume that $d> n \geq 4$ or that $n=3$ and $d\geq 5$. For
$$(F,\sO_F(1))\in \sM_{d,n}(\C),$$
and for any family $f:\sX\to S'$ with $F=f^{-1}(s)$ there exists a neighborhood
$S$ of $s$ and an invertible sheaf $\sL$ with $\sL|_F=\sO_F(1)$ and with $(\sX,\sL)\in \sM_{d,n}(S)$.
\end{lemma}
\begin{proof}
The infinitesimal deformations of $F\subset \BP^n$ are given by
$H^0(F,\sO_F(F))$, whereas those of $F$ are classified by $H^1(F,T_F)$.
Using the exact sequence
$$
0 \>>> T_F \>>> T_{\BP^n}|_F \>>> \sO_F(F) \>>>  0
$$
it is sufficient to show that $H^1(F,T_{\BP^n}|_F)=0$. Since
$$
T_{\BP^n}=\Omega^{n-1}_{\BP^n}\otimes \sO_{\BP^n}(n+1)
$$
this follows from the exact sequence
$$
0 \>>> T_{\BP^n}(-d) \>>> T_{\BP^n} \>>> T_{\BP^n}|_F \>>> 0,
$$
and \ref{1.1}, a) for $2<n$ and $d\neq n+1$, and for $d=n+1$ provided $2\neq n-1$.
\end{proof}
\begin{notations}\label{1.4}
We will also consider the moduli stack $\sM_{d,1}$ of families of $d$ disjoint
points in $\BP^1$. So $\sM_{d,1}(S)$ consists of a $\BP^1$ bundle
$\BP(\sE)$ with a subvariety $\sX\subset \BP(\sE)$ \'etale and finite over $S$
of degree $d$. Again, the family is normalized if $\sO_{\BP(\sE)}(\sX)
=\sO_{\BP(\sE)}(d)$.
\end{notations}
Let us recall next the construction and properties of the Jacobian ring of the hypersurface
$F$, mainly due to Carlson, Green and Griffiths. We follow the presentation given in \cite{Iv1}.

One starts with the commutative diagram of exact sequences
\begin{equation}\label{eq1.1}
\begin{CD}
\noarr 0 \noarr 0 \\
\noarr \V VV \V VV\\
0\>>> \Omega^1_{\BP^n} \>>> {\displaystyle\bigoplus^{n+1}\sO_{\BP^n}(-1)} \>>> \sO_{\BP^n} \>>> 0\\
\noarr \V VV \V VV \V V = V\\
0\>>> \Omega^1_{\BP^n}(\log F) \>>> \sV_F (-1) \>>> \sO_{\BP^n} \>>> 0\\
\noarr \V VV \V VV \\
\noarr \sO_F \> = >> \sO_F,
\end{CD}
\end{equation}
where the upper horizontal sequence is the tautological one, and where $\sV_F(-1)$ is defined
by push out. The second horizontal sequence in (\ref{eq1.1}) splits, and
as explained in \cite{Iv1}, Chapter 2 and Proposition 3.7, one obtains an exact sequence
\begin{equation}\label{eq1.2}
0 \>>> \sO_{\BP^n}(-F) \>>> \bigoplus^{n+1}\sO_{\BP^n}(-1) \>>>
\Omega^1_{\BP^n}(\log F) \>>> 0,
\end{equation}
and its dual
$$
0\>>> \Omega_{\BP^n}^{n-1}(\log F)\otimes \sO_{\BP^n}(n+1-d) \>>>
\bigoplus^{n+1}\sO_{\BP^n}(1)\>>> \sO_{\BP^n}(d)\>>> 0
$$
Let us write for simplicity $\Omega^p_{\BP^n}(\log F)(-\mu)$ instead of
$\Omega^p_{\BP^n}(\log F)\otimes\sO_{\BP^n}(-\mu)$, and
$H^q(\sF)=H^q(\BP^n,\sF)$.

As in \cite{Laz}, Appendix B, one obtains a quasi-isomorphism
between
$$\Omega_{\BP^n}^{n-p}(\log F)(-(p+1)d +n + 1)$$
and the wedge product complex
\begin{equation}\label{eq1.3}
0\>>> \bigoplus^{n+1 \choose p}\sO_{\BP^n}(p(1-d))
\>>> \cdots \>>>\bigoplus^{n+1} \sO_{\BP^n}(1-d) \>\psi >> \sO_{\BP^n}.
\end{equation}
Recall that the Jacobian ideal is defined as
$$
J_\mu ={\rm Im}\{\bigoplus^{n+1}H^0(\sO_{\BP^n}(\mu+1-d))\>\psi >>
H^0(\sO_{\BP^n}(\mu))\},
$$
and $R_\mu=H^0(\BP^n,\sO_{\BP^n}(\mu))/J_\mu.$
The multiplication on polynomials defines a multiplication
$R_\mu\times R_\nu \to R_{\mu+\nu}$, obviously surjective for $\mu, \nu \geq 0$.
For $\sigma=(n+1)(d-2)$
$$
R_\bullet=\bigoplus_{\mu=0}^\sigma R_\mu
$$
is a graded ring, called the Jacobian ring. Macaulay's theorem says:
\begin{equation}\label{macaulay}
R_\sigma=\C \mbox{ \ \ and \ \ }R_\mu\otimes R_{\sigma-\mu} \>>> R_\sigma
\end{equation}
is a perfect pairing (see \cite{CG}). For a real number $a$, we will denote the
integral part by $[a]$ and the roundup by $\roundup{a}=-[-a]$.
\begin{lemma}\label{1.5} \
\begin{enumerate}
\item[(a)] For $\mu \geq p(d-1)-n$ and $r\leq p$ the cohomology of the complex
\begin{multline*}
0\>>> \bigoplus^{{n+1}\choose{p}}H^0(\sO_{\BP^n}(\mu+p(1-d)))\>>>
\bigoplus^{{n+1}\choose{p-1}}H^0(\sO_{\BP^n}(\mu+(p-1)(1-d)))\\ \notag
\>>> \cdots \>>>\bigoplus^{n+1} H^0(\sO_{\BP^n}(\mu+1-d))
\>>> H^0(\sO_{\BP^n}(\mu)) \>>>0
\end{multline*}
at $\displaystyle \bigoplus^{{n+1}\choose{r}}H^0(\sO_{\BP^n}(\mu+r(1-d)))$
is isomorphic to
$$
H^{p-r}(\Omega^{n-p}_{\BP^n}(\log F)(\mu+n+1-d(p+1))).
$$
\item[(b)] For $-d < \mu+n+1-d(p+1)\leq 0,$
hence for
$$p=\roundup{\frac{\mu+n+1-d}{d}}$$
one finds
$$
R_\mu\simeq H^p(\Omega^{n-p}_{\BP^n}(\log F)(\mu+n+1-d(p+1))).
$$
\item[(c)] Writing $T^\ell_{\BP^n}(-\log F)=\wedge^\ell T_{\BP^n}(-\log F)$,
one has
$$
R_{\ell d}\simeq H^\ell(\Omega^{n-l}_{\BP^n}(\log F)(n+1-d))=
H^\ell(T^\ell_{\BP^n}(-\log F)).
$$
\item[(d)] For $\nu=d$ the multiplication $R_\mu\times R_d \to R_{\mu+d}$ in $R_\bullet$ corresponds
under the isomorphisms in b) and c) to the cup-product
\begin{multline*}
H^p(\Omega^{n-p}_{\BP^n}(\log F)(\mu+n+1-d(p+1))) \otimes
H^1(T^1_{\BP^n}(-\log F))\\
\>>>
H^{p+1}(\Omega^{n-p-1}_{\BP^n}(\log F)((\mu+d)+n+1-d(p+2))).
\end{multline*}
\end{enumerate}
\end{lemma}
\begin{proof}
$H^q(\Omega^{n-p}_{\BP^n}(\log F)(\mu+n+1-d(p+1)))$
is isomorphic to the $q$-th hypercohomology of the complex (\ref{1.3}). For all $\ell \leq p$ one has
$\mu + \ell (1-d) > -n$. By \ref{1.1}, iii) and iv), the second spectral sequence of the
hypercohomology of (\ref{eq1.3}) degenerates, and one obtains \ref{1.5}, a).
Part b) and c) are special cases of a). The compatibility of the multiplication and the cup-product
in d) follows from an easy local calculation (see \cite{Iv1}).
\end{proof}
\begin{remark}\label{1.6}
For $p=\roundup{\frac{\mu+n+1-d}{d}}$ write
$\mu+n+1-d=pd - \delta$. If $d$ does not divide $\mu+n+1$
one finds $0 < \delta < d$. Obviously for $\sigma=(n+1)(d-2)$
$$
\sigma-\mu+n+1-d=nd-\mu-n-1 = nd - pd - (d-\delta).
$$
The Macaulay duality $R_\mu\simeq
R_{\sigma-\mu}$ is the Serre duality
$$
H^p(\Omega^{n-p}_{\BP^n}(\log F)(-\delta))\simeq
H^{n-p}(\Omega^{p}_{\BP^n}(\log F)(\delta-d)).
$$
If $\delta=0$, hence $\mu=(p+1)d-n-1$, one finds
$$\roundup{\frac{\sigma-\mu+n+1-d}{d}}=n-1-p.$$
In this case one should identify $R_\mu$ with the primitive
cohomology of $F$
$$
H^p(\Omega^{n-p}_{\BP^n}(\log F))= H^p(\Omega^{n-1-p}_{F})_{\rm prim},
$$
and the Macaulay duality becomes
$$
H^p(\Omega^{n-1-p}_{F})_{\rm prim}\simeq H^{n-1-p}(\Omega^{p}_{F})_{\rm prim}.
$$
\end{remark}
\begin{remark}\label{1.6b}
Lemma \ref{1.5}, b), also allows to calculate the dimension of $R_\mu$. Let us just remark,
that for
$\mu < d-1$ one has $R_\mu=H^0(\sO_{\BP^n}(\mu))$, whereas
$$
\dim(R_{d-1})= \dim(H^0(\sO_{\BP^n}(\mu)))-\dim(J_{d-1})=
\dim(H^0(\sO_{\BP^n}(\mu)))-n-1.
$$
In particular, for $\mu=1, \ldots ,d-1$
the dimension of $R_\mu$ is strictly increasing.
\end{remark}
\begin{lemma}\label{1.7} \
Assume that $d\geq n+1$ and that $\ell \leq n-1$. Then there is an inclusion
$R_{\ell d}\to H^\ell(F,T^\ell_F)$,
and for $\ell < n-1$ both are equal, except possibly for $d=n+1$ and $n=2\ell +1$.
\end{lemma}
\begin{proof}
One has
$$
H^\ell(F,T^\ell_F) = H^\ell(F, \Omega_F^{n-1-\ell}\otimes \omega_F^{-1}),
$$
and that the residue map gives exact sequences
\begin{multline}\label{eq1.7}
H^\ell(\Omega_{\BP^n}^{n-\ell}(n+1-d)) \>>>
H^\ell(\Omega_{\BP^n}^{n-\ell}(\log F)(n+1-d)) \>>>\\
H^\ell(F, \Omega_F^{n-1-\ell}\otimes \omega_F^{-1})\>>>
H^{\ell+1}(\Omega_{\BP^n}^{n-\ell}(n+1-d)).
\end{multline}
By \ref{1.1}, for $\ell \leq n-1$ and for $d>n+1$ the left hand
side in (\ref{eq1.7}) vanishes. For $d=n+1$ the map
$$
H^{\ell-1}(F, \Omega_F^{n-1-\ell})\>>>
H^{\ell}(\Omega_{\BP^n}^{n-\ell})
$$
is surjective.

For $\ell < n-1$ and $d>n+1$ the right hand side of (\ref{eq1.7})
vanishes as well, whereas for $d=n+1$ one needs the additional
condition $\ell+1 \neq n-\ell$.
\end{proof}
The following elementary calculations will be crucial for estimating the length of the Griffiths-Yukawa coupling.
\begin{lemma}\label{1.8} Assume that $d\geq n+1$.
\begin{enumerate}
\item[(a)] The product map
$$
R_\mu\otimes S^{n-1}(R_d) \>>> R_{\mu+(n-1)d}
$$
is non-trivial if and only if $0\leq \mu \leq 2d-2(n+1)$.
\item[(b)] The product map
$$
R_\mu \otimes S^n(R_d) \>>> R_{\mu+nd}
$$
is zero for all $\mu\geq 0$ if and only if $d < 2(n+1)$.
\item[(c)] Let $V$ be a subspace of $R_d$, and let
$$\varsigma(V\subset R_d,\mu)=\varsigma(V,\mu)$$
denote the largest integer $\nu$ for which the product map
$$R_{\mu}\otimes S^\nu(V) \>>> R_{\mu+\nu d}$$
is non-zero (or $0$ if this map is always zero). If $\varsigma(V,d-n-1)<n-1$ then
\begin{gather*}
\varsigma=\varsigma(V,0)= \cdots = \varsigma(V,2(d-n-1))\\
\mbox{and \ } \varsigma \geq \varsigma(V,2(d-n-1)+1) \geq \cdots \geq
\varsigma(V,\sigma)=0.
\end{gather*}
\end{enumerate}
\end{lemma}
\begin{proof} In a) and b) we know, that the product maps are surjective.
Hence a) follows, since $R_{\mu+(n+1)d}=0$, if and only if $\sigma \geq {\mu+(n+1)d}$.

For b) remark that $2(n+1) > d$ is equivalent to $nd > \sigma$.
Hence there exists some $\mu > 0$ with $R_{\mu+nd}\neq 0$ if and only if
$d< 2(n+1)$.

Consider for $\mu$ and $\rho \geq 0$ the commutative diagram
\begin{equation}\label{eq1.8}
\begin{CD}
R_\rho\otimes R_{\mu}\otimes S^\nu(V) \>>> R_\rho \otimes R_{\mu+\nu d}\\
\V\alpha V V \V V V\\
R_{\mu+\rho}\otimes S^\nu(V) \>>> R_{\mu+\nu d +\rho}.
\end{CD}
\end{equation}
Since $\alpha$ is surjective, one finds
\begin{equation}\label{uneq1}
\varsigma(V,\mu) \geq \varsigma(V,\mu+\rho) \mbox{ \ \ for \ \ } \rho \geq 0.
\end{equation}
The multiplication map
$$
R_\mu \otimes V \otimes S^\nu(V) \>>> R_{\mu}\otimes S^{\nu+1}(V)
\>>> R_{\mu+\nu d +d}
$$
factors through the composite
$$
R_\mu \otimes V \otimes S^\nu(V) \>>> R_\mu \otimes R_d \otimes S^\nu(V)
\>>> R_{d+\mu}\otimes S^\nu(V) \>>> R_{\mu+\nu d +d},
$$
hence
\begin{equation}\label{uneq2}
\varsigma(V,\mu)-1 \leq \varsigma(V,\mu+d).
\end{equation}
For $\varsigma=\varsigma(V,0),$ the image $W$ of
$R_0\otimes S^\varsigma(V)= S^\varsigma(V)$ in $R_{\varsigma d}$
is non zero. By the Macaulay duality the map
$R_{\sigma-\varsigma d}\otimes W \>>> R_\sigma$
is surjective, and (\ref{eq1.8}) for $\mu=0$, $\nu=\varsigma$ and
$\rho=\sigma-\varsigma d$ implies that
$\varsigma(V,\sigma - \varsigma d)\geq \varsigma.$ Together with (\ref{uneq1})
one finds that $\varsigma=\varsigma(V,\mu)$ for $\mu=0,\ldots,\sigma - \varsigma d$.

By assumption $\varsigma(V,d-n-1) < n-1$, hence (\ref{uneq2}) implies that
$$
\varsigma \leq \varsigma(V,d)+1 \leq \varsigma(V,d-n-1)+1 \leq n-1,
$$
and $\sigma - \varsigma d \geq \sigma-(n-1)d \geq 2(d-n-1).$
\end{proof}
\section{Cyclic coverings}\label{cyclic}
Let $W$ be a projective manifold, let $\sN$ be an invertible sheaf,
and let $\sigma$ be a section of $\sN^d$ whose zero divisor
$D$ has normal crossings. As in \cite{EV1} one has a normal cyclic covering given by
\begin{equation}\label{sheaves}
Z={\rm \bf Spec}(\bigoplus_{i=0}^{d-1}{\sN^{(i)}}^{-1}) \mbox{ \ \
with \ \ }\sN^{(i)}=\sN^{i}(-[\frac{i\cdot D}{d}]).
\end{equation}
Here $[\frac{i\cdot D)}{d}]=\rounddown{\frac{i\cdot D}{d}}$
denotes the integral part of the $\Q$-divisor $\frac{i\cdot
D}{d}$. We will call $Z$ the variety, obtained by taking the
$d$-th root out of $\sigma$ (or $D$). If $\delta: Z' \to W$ is a
desingularization of $Z$, such that $\delta^*D$ is a normal
crossing divisor, then by \cite{EV1}, 3.22, one has
\begin{equation}\label{sheaves2}
R^b\delta_*\Omega^a_{Z'}(\log (\delta^* D))=\left\{
\begin{array}{ll}\Omega^a_W(\log D)\otimes \bigoplus_{i=o}^{d-1}\sN^{(i)^{-1}}
& \mbox{ \ for \ \ } b=0\\
0 & \mbox{ \ for \ \ } b>0.
\end{array}\right.
\end{equation}
In terms of function fields,
$Z$ is the normalization of $W$ in the Kummer extension
$$\C(W)\big(\sqrt[d]{\frac{D}{d\cdot c_1(\sN)}}\big)$$
where $c_1(\sN)$ stands for a Cartier divisor given by a meromorphic section of
$\sN$ and where $\frac{A}{B}$ denotes a function $f\in \C(W)$
with divisor $A-B$.

Consider for $\iota >0$ the bundle $\pi_\iota:\BP_\iota=\BP(\sO_W\oplus \sN^\iota)\to W.$

For $\iota=d$ there are three natural maps
\begin{gather*}
{\rm id}\oplus \sigma:\sO \>>> \sO\oplus \sN^d, \ \ \ \ \
{\rm id}\oplus 0:\sO \>>> \sO\oplus \sN^d,\\
\mbox{and \ \ }
0\oplus {\rm id}:\sN^d \>>> \sO\oplus \sN^d,
\end{gather*}
inducing sections $s_\sigma$, $s_0$ and $s_\infty$ of $\pi_d:\BP_d \to W$, respectively.
We will write $E_\bullet=s_\bullet(W)$ for the image.

The divisors $E_\infty$ and $E_\sigma$ do not meet,
whereas $E_0\cap E_\sigma$ is isomorphic to the zero divisor $D$ of $\sigma$.
Remark that $E_0+E_\sigma$ is a normal crossing divisor if $D$ is
non-singular. One has
$$
\sO_{\BP_d}(E_\sigma)=\sO_{\BP_d}(E_0)=\sO_{\BP_d}(1), \mbox{ \ \ and \ \ }
\sO_{\BP_d}(E_\infty)=\sO_{\BP_d}(1)\otimes\pi_d^*\sN^{-d}.
$$
The map
$$ \pi_1^*(\sO_W\oplus\sN^d) \>>> S^d(\pi_1^*(\sO_W\oplus\sN))\>>>
\sO_{\BP_1}(d)
$$
defines a morphism $\mu:\BP_1 \to \BP_d$ of degree $d$,
which is the cyclic covering obtained by taking the
$d$-th root out of $(d-1)E_0+E_\infty$ or, using the notation introduced above,
it is the Kummer covering defined by
$$
\sqrt[d]{\frac{(d-1)E_0+E_\infty}{d\cdot c_1(\sO_{\BP_d}(d)\otimes\pi_d^*\sN^{-d})}}=
\sqrt[d]{\frac{E_\infty+d\cdot c_1(\pi_d^*\sN)}{E_0}}.
$$
If $D$ is reduced, hence $\sN^{(i)}=\sN^i$, for $i=0\ldots, d-1$, one finds
$$Z=\mu^{-1}(E_\sigma)\subset \BP_1.$$

Assume that $(f: \sX \to S, \sL)\in \sM_{d,n}(S)$ is a normalized family
of hypersurfaces. For the locally free sheaf $\sE=f_*(\sL)$ on $S$
we choose $W$ to be the total space of the $\BP^n$-bundle
$p:\BP(\sE)\to S$, and $\sN=\sO_{\BP(\sE)}(1)$. The divisor $D=\sX$
given by a section $\sigma$ of $\sO_{\BP(\sE)}(d)$, and taking the $d$-th root
one obtains the cyclic covering $\sZ\to \BP(\sE)$ of degree $d$. As explained
above, $\sZ$ is embedded in $\BP_1$.

The map
$$
\pi_1^*p^*(\sO_S\oplus \sE) \>>> \pi_1^*(\sO_{\BP(\sE)} \oplus \sO_{\BP(\sE)}(1))
\>>> \sO_{\BP_1}(1)
$$
defines a morphism $\eta:\BP_1 \to \BP(\sO_S\oplus \sE)$, which contracts
$E_\infty$ to the section of $\BP(\sO_S\oplus \sE)$ given by
$0\oplus {\rm id}:\sE \to \sO_S\oplus \sE$, and which is an embedding
elsewhere.

In particular, $\eta$ defines an embedding of $\sZ$ into the $\BP^{n+1}$-bundle
$\BP(\sO_S\oplus \sE)$. Since
$$
p_*\pi_{1*}(\sO_{\BP_1}(1)|_Z) = p_* (\sO_{\BP(\sE)}\oplus \sO_{\BP(\sE)}(1))
= \sO_S\oplus \sE
$$
this embedding is defined by the relative sections of $\sO_{\BP_1}(1)|_\sZ$.
Altogether one obtains
\begin{lemma}\label{2.1}
Let $(f: \sX \to S, \sL)\in \sM_{d,n}(S)$ be a normalized
family of hypersurfaces, let $\sE=f_*\sL$ and let
$\pi_1:\sZ \to \BP(\sE)$ be the covering obtained by taking the
$d$-th root out of $\sX$. Then $g=p\circ \pi_1:\sZ\to S$ with the polarization
$\pi_1^*(\sO_{\BP(\sE)}(1))$ is again a normalized family,
embedded in $\BP(\sO_S\oplus \sE)$. The section of $\BP(\sO_S\oplus \sE)$,
given by $0\oplus {\rm id}: \sE \to \sO_S \oplus \sE$,
is disjoint from $\sZ$. Blowing it up one obtains an embedding of
$\sZ$ in the $\BP^1$ bundle $\BP_1=\BP(\sO_{\BP(\sE)}\oplus \sO_{\BP(\sE)}(1))$.
\end{lemma}
Lemma \ref{2.1} allows to iterate the construction of cyclic coverings.
Starting with a normalized family
$$(f: \sZ_0=\sX \to S, \sL)\in \sM_{d,n}(S),$$
and writing $\sE=f_*\sL$ we obtain new families
$$
(g_\nu:\sZ_\nu \to S, \sN_\nu) \in \sM_{d,n+\nu}(S),
$$
by taking successively the $d$-th root out of $\sZ_{\nu-1} \subset
\BP(\sO^{\oplus \nu-1}_S\oplus \sE)$.
\begin{notation}\label{2.2}
We will call $(g_\nu:\sZ_\nu \>>> S, \sN_\nu) \in \sM_{d,n+\nu}(S)$
the family obtained as the $\nu$-th iterated $d$ fold covering of
$$(f: \sZ_0=\sX \to S, \sL)\in \sM_{d,n}(S).$$
Again, we will often write $(g_\nu:\sZ_\nu \to S) \in \sM_{d,n+\nu}(S).$

Let $\sM_{d,n+1}^{(\nu)}$ be the moduli stack of families of
hypersurfaces of degree $d$ in $\BP^{n+\nu}$, obtained in this way, and
let $M^{(\nu)}_{d,n+\nu}$ be the image of $\sM_{d,n+\nu}^{(\nu)}$ in the
coarse moduli scheme $M_{d,n+\nu}$.
\end{notation}
\begin{lemma}\label{2.3}
The morphism $M_{d,n}\to M^{(\nu)}_{d,n+\nu}$ is quasi-finite.
\end{lemma}
\begin{proof}
Let $X\subset \BP^n$ be a smooth hypersurface of degree $d$ and let
$\pi_1:Z_1\to \BP^n$ be the cyclic covering obtained by taking the
$d$-th root out of $X$.

The infinitesimal deformations  of $Z_1$ are given by
$H^1(Z_1,T_{Z_1})$, whereas the ones of the pair $(\BP^n,X)$ are
given by $H^1(\BP^n, T_{\BP^n}(-\log X))$. Since $T_{\BP^n}(-\log
X)$ is a direct factor of $\pi_{1*}T_{Z_1}$ (see \cite{Iv1}, 3.1)
the natural map
$$
\varphi:H^1(\BP^n, T_{\BP^n}(-\log X))\>>> H^1(Z_1,T_{Z_1})
$$
is injective. By \ref{1.7} $H^1(Z_1,T_{Z_1})$ contains
$H^1(\BP^n, T_{\BP^n}(-\log X))$ and by construction
the image of $\varphi$ lies in this subspace. Repeating this argument,
one obtains an injection
$$
H^1(\BP^n, T_{\BP^n}(-\log X))\>>> H^\nu(Z_\nu,T_{Z_\nu}).
$$
In particular the fibres of
$M_{d,n}\to M^{(\nu)}_{d,n+\nu}$ are zero-dimensional.
\end{proof}
Lemma \ref{1.8} allows to give the first estimates for the length
$\varsigma(g)$ of the Griffiths-Yukawa coupling for universal families for $\sM^{(\nu)}_{d,n}$.
\begin{lemma}\label{2.4}
Let $f: \sX \to S$ be a universal family for $\sM_{d,n}$,
and let for $\nu \in \{1,\ldots , d-n-1\}$
$$(g_\nu:\sZ_\nu \to S)\in \sM_{d,n+\nu}$$
be the family obtained as the $\nu$-th iterated $d$ fold covering of $f: \sX \to S$.
\begin{enumerate}
\item[(i)] $\varsigma(f)=n-1$.
\item[(ii)] $\varsigma(g_1)=n-1$, for $n+1 < d < 2(n+1)$.
\item[(iii)] $\varsigma(g_1)=n$ for $d \geq 2(n+1)$.
\item[(iv)] Assume for some $1\leq \nu < d-n$ one has $\varsigma(g_\nu)< n-1+\nu$.
Then $\varsigma(g_\nu) = \varsigma(g_{\nu+1})$.
\end{enumerate}
\end{lemma}
\begin{proof} Let us write $\bigoplus R_\bullet^{(\nu)}$ for the Jacobian ring of
the general fibre $Z_\nu$ of $g_\nu:\sZ_\nu \to S$.

The tangent vectors to $S$ in $s$ are given by $R_d=H^1(\BP^n,T_{\BP^n}(-\log X))$
and i) follows from \ref{1.8}, a), for $\mu=d-n-1$.

As we have seen in the proof of \ref{2.3} $R_d$ is contained in
$$
R^{(1)}_d=H^1(\BP^{n+1},T_{\BP^{n+1}}(-\log Z_1))\subset H^1(Z_1,T_{Z_1}),
$$
hence by induction also in
$$
R^{(\nu-1)}_d=H^1(\BP^{n+\nu-1},T_{\BP^{n+\nu-1}}(-\log Z_{\nu-1}))\subset
H^1(Z_{\nu-1},T_{Z_{\nu-1}}).
$$
The Galois action for the covering $Z_\nu\to \BP^{n+\nu-1}$ induces a decomposition
$$H^{q,p}(Z_\nu)= H^{q,p}(\BP^{n+\nu-1})\oplus \bigoplus_{i=1}^{d-1}
H^{q,p}(Z_\nu)_i$$
in eigenspaces, and the action of $R_d=H^1(\BP^n,T_{\BP^n}(-\log X))$ respects the
decomposition.

For $i\neq 0$ one has (see \cite{EV1}, for example)
$$
H^{q,p}(Z_\nu)_i=H^p(\BP^{n+\nu-1},\Omega_{\BP^{n+\nu-1}}^{q}(\log F)(-i)),
$$
and by \ref{1.5}, b), for $\mu=d(p+1)-i-n-\nu$, hence
for $p=\roundup{\frac{\mu+n+\nu-d}{d}}$,
and for $q=n+\nu-1$ one finds $H^{q,p}(Z_\nu)_i=R^{(\nu-1)}_\mu$. One finds
\begin{equation}\label{eq2.1}
H^{q,p}(Z_\nu)=H^{q,p}(\BP^{n+\nu-1})\oplus \bigoplus_{i=1}^{d-1}
R^{(\nu-1)}_{d(p+1)-i-n-\nu}.
\end{equation}
The cup-product with $R_d \subset H^1(\BP^{n+\nu-1},T_{\BP^{n+\nu-1}}
(-\log Z_{\nu-1}))$ is trivial on the first factor, and it is induced by
the multiplication with $R_d\subset R^{(\nu-1)}_d$ on the others.

Let us consider first the case $\nu=1$ and $d-n-1 \geq 1$. By \ref{1.8},
a), the multiplication $R_\mu\otimes S^{n-1}(R_d) \to
R_{(n-1)d+\mu}$ is nonzero for $\mu=0, \ldots , 2(d-n-1)$, hence
for at least one of the $R_\mu$, occurring in the decomposition (\ref{eq2.1}).
Also, $R_\mu\otimes S^{n}(R_d) \to R_{nd+\mu}$ is zero whenever
$d$ divides $\mu+n+1$, hence for all $R_\mu$ not in (\ref{eq2.1}).
Then ii) and iii), follow directly from \ref{1.8}, a) and b).

Assume now, that for some $d-n-1 \geq \nu >0$ one has
$\varsigma(g_\nu)< n-1+\nu$. For $V=R_d$ as a subspace of $R^{(\nu)}_d$,
using the notation introduced in \ref{1.8}, c), one finds
$$
\varsigma=\varsigma(V\subset R^{(\nu)},d-n-1)=\varsigma(g_\nu)< n-1+\nu=\dim(Z_\nu).
$$
Hence the assumptions made in \ref{1.8}, c), hold true, for $R_\bullet^{(\nu)}$ instead of
$\R_\bullet$, and
$$
R^{(\nu)}_\mu\otimes S^{\varsigma}(R_d) \>>> R^{(\nu)}_{\varsigma d+\mu}
$$
is non zero for $\mu=0, \ldots , 2(d-n-1-\nu)$, hence at least for
one of the $R_\mu$ occurring in the decomposition (\ref{eq2.1}).
On the other hand
$$
R^{(\nu)}_\mu\otimes S^{\varsigma+1}(R_d) \>>> R^{(\nu)}_{\varsigma d+d+\mu}
$$
is always zero. This implies that
$\varsigma(g_{\nu+1})=\varsigma$, as claimed in iv).
\end{proof}
Lemma \ref{2.4} allows to proof the second part of Theorem \ref{6.5}.
\begin{corollary}\label{2.5} For $n+1\leq d < 2n$, and
for $0\leq \ell \leq n-[\frac{d}{2}]$ one has
$$\varsigma(\sM_{d,n}^{(\ell)})=n-\ell-1.$$
\end{corollary}
\begin{proof} Remark that a universal family $g_\ell:\sZ_\ell \to S$
for $\sM^{(\ell)}_{d,n}$ is obtained as the $\ell$-th iterated $d$-fold covering
of $f:\sX \to S\in \sM_{d,j}$, for $[\frac{d}{2}] \leq j=n-\ell \leq n$.
This implies that $2(j+1)>d$, and replacing in Lemma \ref{2.4}, iii) and iv), $n$ by $j$
and one finds $\varsigma(g_1) = \cdots = \varsigma(g_\eta)=j-1$, for
$\eta = 1 ,\ldots , d-j-1$. The condition $d\geq n+1$ implies that
$d-j-1 \geq n-j=\ell$.
\end{proof}

\section{Product subvarieties of the moduli stack}\label{iterated}
\begin{proposition} \label{3.1}
Let $f: \sX_0 \to S=S_1 \times \cdots \times S_m$ be a smooth family
of polarized $m$-folds with general fibre $F$, such that $\dim S_i\geq 1$, for $i=1,\ldots , m$,
and such that the induced map from $S_1 \times \cdots \times S_m$ to the moduli space is
generically finite. Then $\dim S_i\leq h^0(F,\omega^2_F)$, for $i=1,\ldots , m$.

In particular, if $f: \sX_0\to S_1 \times \cdots \times S_m$ is a smooth
family of Calabi-Yau $m$-folds $\dim S_1=\cdots =\dim S_m = 1$.
\end{proposition}
Remark that by \cite{VZ3}, Corollary 6.4, there can not exist families over products with more
than $m$ components.
\begin{proof}[Proof of \ref{3.1}] Let $Y_i$ be a non-singular projective compactification of $S_i$
with $\Gamma_i=Y_i\setminus S_i$ a normal crossing divisor. Choose any extension
$$f: \sX \>>> Y= Y_1\times \cdots \times Y_m$$
of $\sX_0 \to S$ with $\sX$ projective manifold, and with
$\Delta=\sX\setminus \sX_0$ a normal crossing divisor. By construction,
$\Gamma = Y \setminus S$ is a normal crossing divisor, as well.

As in \cite{VZ3}, Section 4, consider the sheaves
$$F^{p,q}=R^qf_*(\Omega^p_{\sX/Y}(\log \Delta)\otimes\omega^{-1}_{\sX/Y}), \quad p+q=m$$
with the Higgs (or Kodaira-Spencer) maps
$$ \tau_{p,q}: F^{p,q}\>>> F^{p-1,q+1}\otimes\Omega^1_Y(\log \Gamma).$$
They induce the $q$-th iterated Kodaira-Spencer map
$$\tau^q=\tau_{m-q+1, q-1}\circ\cdots\circ\tau_{m,0}
: F^{m,0}=\sO_Y\>>> F^{m-q,q}\otimes S^q(\Omega^1_Y(\log \Gamma)).$$
By \cite{VZ3}, Proof of Corollary 6.4, $Y$ can only be the product of $m=\dim F$ non trivial
factors if the $m$-th iterated Kodaira-Spencer map $\tau^m \not=0$.
One obtains an injection
$$ \sO_Y \>>> R^mf_* \omega_{\sX/Y}^{-1}\otimes S^m(\Omega^1_Y(\log \Gamma)),$$
hence, writing $p_i$ for the projection to the $i$-th factor, a non trivial map
\begin{multline}\label{eq3.1}
 \varphi: f_*\omega^2_{\sX/Y} \>>> S^m(\Omega^1_Y(\log \Gamma))=\\
\bigoplus_{j_1+\cdots+j_m=m}S^{j_1}p_1^*\Omega^1_{Y_1}(\log \Gamma_1)
\otimes\cdots\otimes S^{j_m}p_m^*\Omega^1_{Y_m}(\log \Gamma_m).
\end{multline}
By \cite{VZ3} the sheaf $f_*\omega^2_{\sX/Y}$ is big (or ample with respect to some dense open
subset), hence the image $\sK$ of $\varphi$ is big, as well. Since the restriction of each of
the direct factors in (\ref{eq3.1}) to the general fibre of one of the projections is trivial,
except the one for $j_1=\cdots =j_m=1$, the map $\varphi$ factors through a non-trivial map
$$ \varphi: f_*\omega^2_{\sX/Y} \>>> p^*_1\Omega^1_{Y_1}(\log \Gamma_1)
\otimes\cdots\otimes p_m^*\Omega^1_{Y_m}(\log \Gamma_m).$$
For $i=1,\ldots,m$, the restriction of the latter to a general fibre of the $i$-th projection
gives a non-trivial map
$$ f_*\omega^2_{\sX/Y}|_{Y_j}\>>> \bigoplus \Omega^1_{Y_j}(\log \Gamma_j),$$
and the composite with the projection to one of the direct factors must be non-trivial.
Since $f_*\omega^2_{\sX/Y}|_{Y_j}$ is again big, Bogomolov's lemma implies that
this projection is surjective on some open dense subset. In particular,
$$ {\rm rank} (\Omega^1_{Y_j}(\log \Gamma_j)) \leq
{\rm rank}(f_*\omega^2_{\sX/Y}|_{Y_j})=h^0(F,\omega^2_F).$$
\end{proof}\ \\
The bound in \ref{3.1} is far from being optimal. In fact, one could hope that
one has the following injectivity of wedge products of tangent
spaces.
\begin{conjecture}\label{3.2}
Let $f: \sX_0 \to S=S_1 \times \cdots \times S_\ell=S$ be a smooth family
of polarized manifolds  with general fibre $F$, such that $\dim S_i\geq 1$, for $i=1,\ldots , \ell$,
and such that the induced map from $S_1 \times \cdots \times S_\ell$ to the moduli space is
generically finite. Then for all $1\leq  k\leq \ell$ the composition
$$ \bigoplus_{1\leq i_1<\cdots< i_k\leq \ell}
p^*_{i_1}(T_{S_{i_1}})\otimes\cdots\otimes p^*_{i_k}(T_{S_{i_k}})\>>>\bigoplus^{k \choose \ell}
(R^1f_*T_{\sX_0/S})^{\otimes k}\>>> R^k f_*T^k_{\sX_0/S}$$
is injective. In particular
$$ \sum_{1\leq i_1<\cdots<i_k\leq \ell}\dim(S_{i_1})\cdot \cdots \cdot \dim(S_{i_k}) \leq h^k(F,T^k_F)
$$
for a general fibre $F$ of $f$ and for $T^k_\bullet=\wedge^k T_\bullet$.
\end{conjecture}
Note that over a local bases one can easily find examples, where for $k>0$ the above
wedge product map is not injective.

If the general fibre of $f$ is a Calabi-Yau $m$-fold, or more generally a $m$-dimensional
manifold of Kodaira dimension zero, then the conjecture follows from \ref{3.1}, for
$k=\ell=m$. In fact, as we have seen in the proof of \ref{3.1} the $m$-th iterated Kodaira-Spencer
map has to be non-trivial, and its image must lie in
$$
p^*_1\Omega^1_{Y_1}(\log \Gamma_1)
\otimes\cdots\otimes p_m^*\Omega^1_{Y_m}(\log \Gamma_m).
$$
So the map in \ref{3.2} is non trivial, and by \ref{3.1} it is injective.

We will give below an affirmative answer to \ref{3.2}
for families of Calabi-Yau manifolds, and for families  of hypersurfaces in $\BP^n$
of degree larger than or equal to $n+1$.\\

Let $\V$ be an irreducible  polarized complex variation of Hodge structures over the
product of quasi-projective manifolds $S_1\times \cdots \times S_\ell$, and let
$$p_i: S_1\times \cdots\times S_\ell\>>> S_i$$
denote the projections.
\begin{proposition}\label{3.3}
There exist polarized complex variation of Hodge structures  $\V_i$ on $S_i,\, 1\leq i\leq \ell$,
and a Hodge isometry $ \V\simeq p_1^*(\V_1)\otimes \cdots\otimes p^*_\ell(\V_\ell).$
\end{proposition}
\begin{proof} We have to consider the case $\ell=2.$ The general case
follows by induction.
Let $\rho: \pi_1(S_1)\times \pi_1(S_2)\to {\rm Gl}_n(\C)$ denote
the underlying representation of the fundamental group. Let $G_1$ and $G_2$
be the Zariski closure of the image of the representations
$$\rho: \pi_1(S_1)\times\{e_2\}\>>> {\rm Gl}_n(\C),\mbox{ \ \ and \ \ }
\rho: \{e_1\}\times\pi_1(S_2)\>>> {\rm Gl}_n(\C),$$
respectively. By \cite{De1} $G_1$ and $G_2$ are semi simple algebraic groups.
By construction  $\rho$ factors through the natural representation of
$G_1\times G_2$ in ${\rm Gl}_n(\C),$
$$
\begin{CDS}
\pi_1(S_1)\times\pi_1(S_2) \>  \rho   >> {\rm Gl}_n(\C)\\
\novarr \SE E (\tau_1,\tau_2) E  \quad\A A \gamma A  \\
\novarr \novarr \quad
G_1\times G_2.\end{CDS}
$$
By Schur's Lemma the representation $\gamma$ can be decomposed as tensor product of
representations $\gamma=\gamma_1\otimes \gamma_2$
with 
$$\gamma_1:G_1\>>> {\rm Gl}_{n_1}(\C),\mbox{ \ \  and \ \ }\gamma_2:G_2\>>> {\rm Gl}_{n_2}(\C).$$
In this way one obtains a decomposition $\rho=\rho_1\otimes\rho_2$ with
$$\rho_1: \pi_1(S_1)\>>> {\rm Gl}_{n_1}(\C),\mbox{ \ \ and \ \ }\rho_2: \pi_1(S_2)\>>> {\rm Gl}_{n_2}(\C).$$
\begin{claim}\label{3.4}
The local systems $\V_1$ and $\V_2$, given by $\rho_1$ and $\rho_2$, respectively, are underlying
polarized complex variations of Hodge structures.
\end{claim}
\begin{proof}
By \cite{De1}, p. 9, the local system given by the restriction of $\rho$ to $\pi_1(S_1)\times
\{e_2\}$ decomposes as a direct
sum of irreducibles polarized complex variation of Hodge structures. Each of them
corresponds to an irreducible representation $\rho^j$ over $\C$. Since
$$ \rho_1^{\oplus n_2}=\rho|_{\pi_1(S_1)\times\{e_2\}}=\sum_j\rho^j$$
one finds $\rho^j=\rho_1$ for all $j.$ In particular, $\V_1$ admits a polarized
complex variation of Hodge structures.
\end{proof}
By Claim \ref{3.4} and by the functoriality of polarized complex variation of Hodge structures,
$\V_1\otimes\V_2$ is a local system, underlying a polarized complex variation of Hodge structures.
Note that $\rho=\rho_1\otimes\rho_2$ is irreducible over $\C$. By \cite{De1}, Prop. 1.13,
(or by \cite{Sim1}, Lemma 4.1, in the projective case) a $\C$-irreducible
local system admits at most one polarized complex variation of Hodge structures. Hence the
original polarized complex variation of Hodge structures on $\V$ coincides with
the one coming from the tensor product $\V_1\otimes\V_2.$
\end{proof}

\begin{proposition} \label{3.5} $f: \sX\to S=S_1\times \cdots\times S_\ell$ be a smooth
family of minimal $m$-folds, such that the induced map from
$S_1\times\cdots\times S_\ell$ to the corresponding moduli space is generically finite.
\begin{enumerate}
\item[(a)] If the fibres of $f$ are Calabi-Yau $m$-folds (or more generally,
if $\Omega^m_{\sX/S}$ is the pullback of an invertible sheaf on $S$),
the composition
\begin{multline*}
\bigoplus_{1\leq i_1<\cdots< i_k\leq \ell}
p^*_{i_1}(T_{S_{i_1}})\otimes\cdots\otimes p^*_{i_k}(T_{S_{i_k}})\>>>\bigoplus^{k\choose \ell}
(R^1f_*T_{\sX_0/S})^{\otimes k}\>>> R^k f_*T^k_{\sX_0/S}
\end{multline*}
is injective for $1\leq k\leq \ell.$ In particular,
$$ \sum_{1\leq i_1<\cdots<i_k\leq \ell}\dim(S_{i_1})\cdot \cdots \cdot \dim(S_{i_k})
\leq h^k(F,T^k_F)=h^k(F,\Omega^{m-k}).$$
\item[(b)] Assume that $f: \sX\to S=S_1\times \cdots\times S_\ell$ is a smooth normalized
family of hypersurfaces in $\BP^{n}$ of degree $d \geq n+1$.
Then the natural maps
\begin{multline*}
\bigoplus_{1\leq i_1<\cdots< i_k\leq \ell}p^*_{i_1}(T_{S_{i_1}})\otimes\cdots\otimes
p^*_{i_k}(T_{S_{i_k}}) \>>>\\
 R^k p_*(\Omega^{n-k}_{\BP(\sE)/S}(\log \sX)(n+1-d))\otimes
\det(\sE)^{-1} \>>> R^k f_*T^k_{\sX_0/S}
\end{multline*}
are both injective for all $1\leq k\leq \ell.$ In particular,
\begin{multline*}
\sum_{1\leq i_1<\cdots<i_k\leq \ell}\dim(S_{i_1})\cdot \cdots \cdot \dim(S_{i_k})
\leq\\
h^k(\BP^n,\Omega^{n-k}_{\BP^n}(\log F)(-(d-(n+1)))
\leq h^k(F,T_F^k).
\end{multline*}
\end{enumerate}
\end{proposition}

\begin{proof}
a) Consider the polarized $\Q$-variation of Hodge structures
$$\V=R^mf_*\Q_{\sX/S}.$$
By assumption, the rank of $f_*\Omega^m_{\sX/S}$
is one, and there exists an irreducible direct factor $\V'$ of $\V$ whose system of Hodge
bundles contains $f_*\Omega^m_{\sX/S}$. Then
$\V'$ is a polarized complex variation of Hodge structures
with the Hodge decomposition
$$\bigoplus_{p+q=m}{E^{p,q}}$$
with $E^{m,0}=f_*\Omega^m_{\sX/S}.$ The Kodaira-Spencer map, injective by assumption,
factors through
$$ T_S\otimes E^{m,0}\>>> {E^{m-1,1}}.$$
By Proposition \ref{3.3} one finds a decomposition of polarized
complex variations of Hodge structures
$$ \V'=p^*_1\V_1\otimes \cdots\otimes p^*_\ell\V_\ell.$$
Let us write the Hodge bundles of $\V_i$ for $i=1,\ldots,\ell$ as
$$ \bigoplus_{p+q=m_i}F^{p,q}_i.$$
Comparing the Hodge bundles in the tensor product decomposition one finds
\begin{gather*}
\bigotimes_{j=1}^\ell p^*_{j}F_{j}^{m_{j},0}=E^{m,0},
\mbox{ \ \ and}\\
\bigoplus_{1\leq i_1<\cdots< i_k\leq \ell}(\bigotimes_{j=1}^k p^*_{i_j}F_{i_j}^{m_{i_j}-1,1}
\otimes \bigotimes_{j=k+1}^\ell p^*_{i_j}F^{m_{i_j},0})\subset E^{m-k,k}.
\end{gather*}
Here we write $\{1, \ldots , \ell\}$ as the disjoint union of $\{i_1, \ldots , i_k\}$
and $\{i_{k+1},\ldots , i_\ell\}$. Then
\begin{gather*}
{\rm rank }F_{1}^{m_1,0}=\cdots={\rm rank} F_{\ell}^{m_\ell,0}=1\mbox{ \ \ and}\\
\sum_{1\leq i_1<\cdots< i_k\leq \ell}{\rm rank}F_{i_1}^{m_{i_1}-1,1}
\cdot \cdots \cdot {\rm rank}F_{i_k}^{m_{i_k}-1,1}\leq {\rm rank} E^{m-k,k}.
\end{gather*}
The injectivity of the Kodaira-Spencer map implies that $F^{m_{i},0}\otimes T_{S_{i}}\to
F^{m_{i}-1,1}$ is injective for $1\leq i\leq \ell,$ hence that
\begin{multline*}
\bigoplus_{1\leq i_1<\cdots< i_k\leq \ell}
p_{i_1}^*(T_{S_{i_1}})\otimes \cdots \otimes p_{i_k}^*(T_{S_{i_k}})
\>>> E^{m-k,k}\otimes (E^{m,0})^{-1}\\
\subset R^kf_*(\Omega_{\sX/S}^{m-k})\otimes
(E^{m,0})^{-1} =R^kf_*T^k_{\sX/S}
\end{multline*}
is injective, as well. One obtains \ref{3.5}, a).\\[.2cm]
b) Consider as in \ref{2.1} for $\sX \subset \BP(\sE)$ the
$d$ fold cyclic cover of $\BP(\sE)$ obtained by taking the $d$-th root out of
$\sX\subset \BP(\sE)$. It gives rise to a family $g: \sZ\to S$ of hypersurfaces in
$\BP(\sO_S\oplus \sE)$ of degree $d$. The Galois action on $\sZ$ induces a decomposition
of $R^ng_*(\C_\sZ)$ into eigenspaces
$$ R^ng_*(\C_\sZ)=\bigoplus_{i=0}^{d-1}R^ng_*(\C_\sZ)_i.$$
For $i\neq 0$, each factor is a polarized complex variation of Hodge structures, with system
of Hodge bundles
$$E^{n-q, q}_i=R^qp_*(\Omega^{n-q}_{\BP(\sE)/S}(\log \sX)(-i)).$$
Observe that for those $i$
$$
E^{n,0}_i\simeq p_* \omega_{\BP(\sE)/S}(d-i) \simeq
p_*\sO_{BP(\sE)}(d-i-(n+1)) = \det(\sE).$$ So for $d>n+1$ the
polarized complex variation of Hodge structures
$$\V= R^ng_*(\C_\sZ)_{d-(n+1)}$$
has the property that $E^{n,0}_{d-(n+1)}\simeq \det(\sE)$.

For $d=n+1$ we consider as in \ref{1.6} the primitive part
$\V$ of the variation of Hodge structures $R^{n-1}f_*(\C_\sX)$,
which is isomorphic to the variation of mixed Hodge structures
$R^np_*(\C_{\BP(\sE)\setminus \sX})$. Again the first
Hodge bundle $E^{n,0}$ is of rank one and isomorphic to $\det(\sE)$.

\begin{claim}\label{3.6}
For $d>n+1$ the Kodaira-Spencer map
$$\theta_{n,0}: E^{n,0}_{d-(n+1)}\otimes T_S \>>> E^{n-1,1}_{d-(n+1)}$$
is injective.
\end{claim}
\begin{proof}
Taking residues along the divisor $\sX\subset\BP(\sE)$ of
\begin{multline*}
\theta^{n,0}: p_*\Omega^n_{\BP(\sE)/S}(\log \sX)(n+1-d)\otimes T_S
\>>>  R^1p_*\Omega^{n-1}_{\BP(\sE)/S}(\log \sX)(n+1-d)
\end{multline*}
one obtains the Kodaira-Spencer map
$$ \tau^{n-1,0}: f_*(\Omega^{n-1}_{\sX/S}\otimes\omega^{-1}_{\sX/S})
\otimes T_S=T_S\>>> R^1f_*(\Omega^{n-2}_{\sX/S}
\otimes\omega^{-1}_{\sX/S}),$$
which is injective. On the other hand, by \ref{1.1} the residue map on
$$ p_*\Omega^n_{\BP(\sE)/S}(\log \sX)(n+1-d)$$
is injective, as well.
\end{proof}
Let $\V'$ be the irreducible direct factor of $\V$ containing $E^{n,0}_{d-(n+1)}$.
Remark that for $d=n+1$ one has a mixed variation of
Hodge structures, but it is isomorphic a polarized variation of Hodge structures
of weight $n-1$.
Hence in both cases by Proposition \ref{3.3} one finds a decomposition of
polarized complex variation of Hodge structures
$$ \V'=p^*_1\V_1\otimes \cdots\otimes p^*_\ell\V_\ell.$$
Let us write again the Hodge bundles of $\V_i$ for $i=1,\ldots,\ell$ as
$$ \bigoplus_{p+q=m_i}F^{p,q}_i.$$
Comparing the Hodge bundles in the above tensor product decomposition one finds
as above
\begin{gather*}
\bigotimes_{j=1}^\ell p^*_{j}F^{m_{j},0}=\det(\sE)
\mbox{ \ and \ }\\
\bigoplus_{1\leq i_1<\cdots< i_k\leq \ell}(\bigotimes_{j=1}^k p^*_{i_j}F_{i_j}^{m_{i_j}-1,1}
\otimes \bigotimes_{j=k+1}^\ell p^*_{i_j}F_{i_j}^{m_{i_j},0})\subset E^{n-k,k}_{d-(n+1)}.
\end{gather*}
Hence,
\begin{gather*}
{\rm rank }F_1^{m_1,0}=\cdots={\rm rank} F_\ell^{m_\ell,0}=1\mbox{ \ \ and}\\
\sum_{1\leq i_1<\cdots< i_k\leq \ell}{\rm rank}F_{i_1}^{m_{i_1}-1,1}
\cdot \cdots \cdot {\rm rank}F_{i_k}^{m_{i_k}-1,1}
\leq {\rm rank} E^{n-k,k}_{d-(n+1)}.
\end{gather*}
Claim \ref{3.6} implies that
$ F^{m_i,0}\otimes T_{S_i} \to F^{m_i-1,1}$ is injective
for $1\leq i\leq \ell$, hence that
$$
\bigoplus_{1\leq i_1<\cdots< i_k\leq \ell}
p_{i_1}^*(T_{S_{i_1}})\otimes \cdots \otimes p_{i_k}^*(T_{S_{i_k}})
\>>> E^{n-k,k}_{d-(n+1)}\otimes \det(\sE)^{-1}
$$
is injective, as well. The sheaf $E^{n-\ell,\ell}_{d-(n+1)}$ is a sub sheaf of
$$R^\ell p_*(\Omega^{n-\ell}_{\BP(\sE)/S}(\log \sX)(n+1-d)),$$
and one obtains the first inclusion in \ref{3.5}.
\begin{claim}\label{3.7}
The residue map
\begin{multline}\label{eq3.2}
R^\ell p_*(\Omega^{n-\ell}_{\BP(\sE)/S}(\log \sX)(n+1-d))\>>>\\
R^\ell f_* (\Omega^{n-\ell-1}_{\sX/S} \otimes\sO_{\sX}(n+1-d))\simeq
R^\ell f_*T_{\sX/S}^\ell\otimes \det(\sE)
\end{multline}
is injective.
\end{claim}
\noindent {\it Proof.} \,
For the isomorphism in (\ref{eq3.2}) remark that the sheaf on the right hand
side is
$$
R^\ell f_* (T^\ell_{\sX/S}\otimes\omega_{\sX/S}(n+1-d))=
R^\ell f_*T_{\sX/S}^\ell\otimes \det(\sE).
$$
For the injectivity it is sufficient to consider a point $s$ in $S$,
and to apply \ref{1.7}.
\end{proof}
\section{Irreducibility of certain local systems and rigidity}\label{rigidity}

Recall that Deligne has shown in \cite{De3} that for the primitive cohomology of
a universal family $f:\sX \to S$ for $\sM_{d,n}$ the monodromy representation is irreducible.
This was extended (although not stated explicitly) by Carlson and Toledo in
\cite{CT} to the local systems $(R^ng_*\C_\sZ)_i$, obtained as eigenspaces
of the variation of Hodge structures of the $d$-fold cyclic covering
$g:\sZ \to S$, for certain values of $d$ and $i$. We use a different argument
to obtain this result for $d < 2(n+1)$. It will be needed to prove Theorem \ref{8.3}.

\begin{lemma}\label{8.0} Assume that $n+1\leq d < 2(n+1)$.
\begin{enumerate}
\item[(i)] For
$$i= 1,\ldots , d-n-1 \mbox{ \ \ and \ \ } i=n+1,\ldots, d-1$$
the polarized complex variations of Hodge structures
$(R^ng_{*}\C_{\sZ/S})_i$ are irreducible.
\item[(ii)] In i) the length of the Griffiths-Yukawa coupling of $(R^ng_{*}\C_{\sZ/S})_i$ is $n-1$,
whereas for $d-n-1 < i < n+1$ it is strictly smaller than $n-1$.
\item[(iii)] For $i\neq 0$ the first Hodge bundle $E_i^{n,0}$ of $(R^ng_{*}\C_{\sZ/S})_i$ is
non-zero, if and only if $i \leq d-n-1$.
\item[(iv)] ${\rm rank}E^{n,0}_{d-n-1} < {\rm rank}E^{n,0}_{d-n-2} < \cdots <
{\rm rank}E^{n,0}_{1}$\\
\hspace*{3cm} $<{\rm rank}E^{n-1,1}_{d-1}< {\rm rank}E^{n-1,1}_{d-2} < \cdots <
{\rm rank}E^{n-1,1}_{n+1} $.
\end{enumerate}
\end{lemma}
\begin{proof}
Remark that \ref{1.5}, b), allows to identify the rank of the sheaves in iv) with
$$
\dim R_0, \ \dim R_1, \ \cdots, \ \dim R_{d-n-2}, \ \dim R_{d-n}, \ \dim R_{d-n+1}, \
\cdots, \ \dim R_{2d-2(n+1)}.$$
Since ${2d-2(n+1)} < d$, iv) follows from \ref{1.6b}.
Part iii) follows from \ref{1.5}, b), and ii) follows from \ref{1.8}, a) and b).

Let us assume that for some $1\leq i \leq d-n-1$ the local system $(R^ng_{*}\C_{\sZ/S})_i$
is a direct sum of irreducible non-trivial sub-systems $\V^1_i\oplus\cdots\oplus\V^\ell_i$, with $\ell>1$.
Using complex conjugation, this is equivalent to the existence of a direct sum decomposition
$\V^1_{d-i}\oplus\cdots\oplus\V^\ell_{d-i}$
of $(R^ng_{*}\C_{\sZ/S})_{d-i}$ in irreducible non-trivial local sub-systems.
Remark that $(R^ng_{*}\C_{\sZ/S})_i\oplus (R^ng_{*}\C_{\sZ/S})_{d-i}$ is defined over $\R$
and polarized by restricting the polarization of $(R^ng_{*}\C_{\sZ/S})_{\rm prim}$. By \cite{De1},
a polarized complex variation of Hodge structures over a quasi-projective manifold
can be written as a direct sum of irreducible ones, orthogonal to each other. In particular,
we may assume that the direct factors $\V^\iota_i\oplus \V^\iota_{d-i}$ are orthogonal with
respect to the polarization. Let us write
$$
\bigoplus_{p=0}^{n-1} E_i^{\iota \, n-p,p}\mbox{ \ \ and \ \ }\bigoplus_{p=1}^n
E_{d-i}^{\iota \, n-p,p}
$$
for the system of Hodge bundles for $\V^\iota_i$ and $\V^\iota_{d-i}$. Let us use again
the abbreviation $\sF(\nu)=\sf\otimes \sO_{\BP(\sE)}(\nu)$. For $j=i$ or $j=d-i$ one has
$$E_j^{\iota \, n-p,p}\subset
R^pp_*(\Omega_{\BP(\sE)/S}^{n-p}(\log \sX)(-j)).$$

Part i) of \ref{8.0} for $d=n+2$, or more general
for $i=d-n-1$, follows from the next Claim, as the $(n,0)$-Hodge bundle of
$(R^ng_{*}\C_{\sZ/S})_{d-n-1}$ is of rank one.
\begin{claim}\label{i1}
For $\iota=1,\ldots,\ell$ one has $E_i^{\iota \, n,0}\neq 0$. In particular
$i\leq d-n-2$.
\end{claim}
\begin{proof}
The restriction of $\bigoplus_{\iota=1}^\ell E_i^{\iota \, n-p,p}$ to a point $s\in S$ is
$R_{(p+1)d-n-1-i}$. Assume that $E_i^{1 \, n,0}=0$, but $E_i^{1 \, n-p,p}\neq 0$, for some $p>0$.
By \ref{1.5}, d), the multiplication with $R_d$ respects the decomposition, and the image of
$$R_{d-n-1-i}\times S^p(R_d) \>>> R_{(p+1)d-n-1-i}$$
lies in
$$
\bigoplus_{\iota=2}^\ell ( E_i^{\iota \, n-p,p})_s,
$$
contradicting the surjectivity of the multiplications maps in the Jacobian ring.
\end{proof}
The cup product on $\V^\iota_i$ induces
$$\Phi_1^\iota:E_i^{\iota \, n-p,p}\otimes R^1p_*(T_{\BP(\sE)/S}(-\log \sX)(2i-d)) \>>>
R^{p+1}p_*(\Omega_{\BP(\sE)/S}^{n-p-1}(\log \sX)(i-d)),
$$
and
$$\Phi_2^\iota:E_i^{\iota \, n-p,p}\otimes p_*(\sO_{\BP(\sE)}(i)) \>>>
R^{p}p_*(\Omega_{\BP(\sE)/S}^{n-p}(\log \sX)).
$$
\begin{claim}\label{8.01}
The image of $\Phi_1^\iota$ is $E_{d-i}^\iota \,{}^{n-p-1,p+1}$.
\end{claim}
\begin{proof}
Remark, that by \ref{1.5} and \ref{1.6}, for all points $s\in S$, the
restriction of the map $\Phi_1^\iota$ to $s$ is induced by the multiplication map
$$
R_{(p+1)d-n-1-i}\otimes R_{2i} \>>> R_{(p+1)d-n-1+i}.
$$
Let us write $V^\iota_{(p+1) d-n-1-i}$ for the subspace of $R_{(p+1)d-n-1-i}$, corresponding
to ${E_i^\iota \,}^{n-p,p}$ and similarly $V^\iota_{p d-n-1+i}$ for the subspace of $R_{pd-n-1+i}$,
corresponding to $E_{d-i}^{\iota \, n-p,p}$.

Since $2(d-n-1-i)\leq d-2i$, one finds by \ref{1.6b} that
$$R_{2(d-n-1-i)}=H^0(\sO_{\BP^n}(2(d-n-1-i))).$$
In particular  $r^2\in R_{2(d-n-1-i)}$ is non-zero, for
$r\in R_{d-n-1-i}\setminus\{0\}$. Assume that for some $t\in R_{2i}$ one has
$t\cdot r\neq 0$. The Macaulay duality, and the surjectivity of the multiplication maps in
$R_\bullet$ allows to choose some $t'\in R_{(n-1)d}$ with $r^2\cdot t'\cdot t\neq 0$ in $R_{\sigma}=
R_{(n+1)(d-2)}$.
Then $(t'\cdot r)\cdot (t\cdot r)\neq 0$. Since $t'\cdot r \in V^\iota_{nd-n-1-i}$ the orthogonality
of the decomposition implies that $t\cdot r\in V^\iota_{d-n-1+i}$.

So the image of $\Phi_1^\iota$ is contained in $E_{d-i}^{\iota \, n-p-1,p+1}$, for all $\iota$.
On the other hand, the multiplication maps $R_\nu\otimes R_\mu$ are surjective for
all $\nu,\mu \geq 0$. Hence ${\rm Im}(\Phi_1^1) \cup \cdots  \cup {\rm Im}(\Phi_1^\ell)$
spans
$$
\bigoplus_{\iota=1}^\ell E_{d-i}^{\iota \, n-p-1,p+1},
$$
and one obtains Claim \ref{8.01}.
\end{proof}
Next define $F^{\iota \, n-p,p}=\Phi_2^\iota(E_i^{\iota \, n-p,p}\otimes p_*(\sO_{\BP(\sE)}(i)))$.
\begin{claim}\label{8.02} For $p=0,\ldots,n-1$
$$
\bigoplus_{\iota=0}^{\ell}F^{\iota \, n-p,p} = R^{p}p_*(\Omega_{\BP(\sE)/S}^{n-p}(\log \sX)).
$$
\end{claim}
\begin{proof}
Again we take a point $s\in S$. By \ref{1.5} and \ref{1.6} the restriction of
$$
\bigoplus_{p=0}^{n-1}R^{p}p_*(\Omega_{\BP(\sE)/S}^{n-p}(\log \sX)
$$
to $s$ is isomorphic to
$$
R_{d-n-1}\oplus R_{2d-n-1} \oplus \cdots \oplus R_{nd-n-1},
$$
whereas the systems of Hodge bundles for $(R^ng_*\C_\sZ)_i$ and $(R^ng_*\C_\sZ)_{d-i}$, restricted
to $s$ are
\begin{gather*}
R_{d-n-1-i}\oplus R_{2d-n-1-i} \oplus \cdots \oplus R_{nd-n-1-i}, \mbox{ \ \ and}\\
R_{d-n-1+i}\oplus R_{2d-n-1+i} \oplus \cdots \oplus R_{nd-n-1+i}.
\end{gather*}
Here we use again $d < 2(n+1)$, or equivalently $\sigma=(n+1)(d-2) < nd$.

For $p=0, \ldots, n-1$ the decompositions of the last two local systems
induce the decompositions
$$
R_{(p+1)d-n-1-i} = \bigoplus_{\iota=1}^\ell V^{\iota}_{(p+1)d-n-1-i}\mbox{ \ \ and \ \ }
R_{(p+1)d-n-1+i} = \bigoplus_{\iota=1}^\ell V^{\iota}_{(p+1)d-n-1+i}
$$
and the image of the multiplication map
$$V^\iota_{(p+1)d-n-1-i}\otimes R_d \>>> R_{(p+2)d-n-1-i}$$
is contained in $V^\iota_{(p+2)d-n-1-i}$. By \ref{8.01} we also know, that the image of
$$V^\iota_{(p+1)d-n-1-i}\otimes R_{2i}\>>> R_{(p+2)d-n-1+i}$$
is contained and generically isomorphic to $V^\iota_{(p+1)d-n-1+i}.$

By our choice of $F^{\iota \, n-p,p}$, its restriction to $s$ is
the image $W^\iota_{(p+1)d-n-1}$ of the multiplication map
$$
V^\iota_{(p+1)d-n-1-i}\otimes R_{i} \>>> R_{(p+1)d-n-1}.
$$
The surjectivity of this map implies that
$$R_{(p+1)d-n-1}= \sum_{\iota=1}^\ell W^\iota_{(p+1)d-n-1}.
$$
For $U_\eta=W^1_{(p+1)d-n-1}\cap \sum_{\iota=2}^\ell W^\iota_{(p+1)d-n-1}$
the image of
$$
U_\eta\otimes R_{i} \>>> R_{(p+1)d-n+i}
$$
lies in $V^{1}_{(p+1)d-n-1+i}\cap \bigoplus_{\iota=2}^\ell V^\iota_{(p+1)d-n-1+i}$,
hence it is zero. Then
$$
U_\eta \times R_{\sigma - (p+1)d+n+1}\>>> R_\sigma
$$
is zero. The Macaulay duality (\ref{macaulay}) implies that $U_\eta=0$, and applying this
to all intersections one finds
$$
R_{(p+1)d-n-1}= \bigoplus_{\iota=1}^\ell W^\iota_{(p+1)d-n-1}.
$$
\end{proof}
As in Remark \ref{1.6} for $E^{n-p,p}=R^pp_*(\Omega_{\BP(\sE)/S}^{n-p}(\log \sX))$ the Higgs bundle
$$
(E,\theta)=\big(\bigoplus_{p=0}^{n-1}E^{n-p,p},\bigoplus_{p=0}^{n-1}\theta_{n-p,p}\big)
$$
is, up to a shift by $(-1,0)$ in the bidegrees, the one corresponding to
$(R^{n-1}f_*\C_\sX)_{\rm prim}$.
By definition of $F^{\iota \, n-p,p}$ the image of $F^{\iota \, n-p,p}\otimes
R^1p_*T_{\sX}$ lies in $F^{\iota \, n-p-1,p+1}$, hence we constructed a decomposition
in a direct sum of sub Higgs bundles
\begin{equation}\label{eq8.01}
\big(\bigoplus_{p=0}^{n-1}E^{n-p,p},\bigoplus_{p=0}^{n-1}\theta_{n-p,p}\big)=
\bigoplus_{\iota=1}^\ell\big(\bigoplus_{p=0}^{n-1}F^\iota \,{}^{n-p,p},\bigoplus_{p=0}^{n-1}
\theta^\iota_{n-p,p}\big).
\end{equation}

Consider next a  general Lefschetz pencil of hypersurfaces of degree $d$
$$\Psi:\tilde \BP^n\subset \BP^n\times \BP^1 \>>> \BP^1$$
smooth over $\BP^1\setminus\{s_1,\ldots,s_r\}$.

By \cite{DK} the fibres of $\Psi$ have at most one rational double point of type $A_1$, i.e.
a singularity given locally analytic as the zero set of the equation $x_1^2+\ldots+x_n^2$ in $\C^n$.

Let $\tau:C\to \BP^1$ be a finite covering, $C_0=\tau^{-1}(\BP^1\setminus\{s_1,\ldots,s_r\})$,
such that the morphism $C_0\to M_{d,n}$ lifts to a morphism $\rho:C_0\to S$. The pullback of the
families over $S$ and over $\BP^1\setminus\{s_1,\ldots,s_r\}$ to $C_0$ are isomorphic.

Replacing $C$ by some covering, one can write $\BP^n\times C = \BP(\sE_C)$ for a locally free
sheaf $\sE_C$ on $C$, such that the closure $\bar\sX_C$ of $\sX\times_SC_0$ is the zero-set of
a section of $\sO_{\BP(\sE_C)}(d)$. Finally we choose an embedded desingularization $\delta_C:\BP_C
\to \BP(\sE_C)$ of $\bar\sX_C$.

Let $D\to \BP^1$ be a $2$-fold covering, totally ramified in the points  $s_1,\ldots,s_r$.
Then $\bar\sX_D$ is locally given by the equation $t^2-x_1^2+\ldots+x_n^2$ in $\C^{n+1}$.
Blowing up the singular point, one obtains an embedded desingularization of $\bar\sX_D$, such that
the total transform $\sX_D+E$ is a normal crossing divisor, $(E)_{\rm red}\cong\BP^n$, and
$E=2E_{\rm red}$. We will assume that $C\to\BP^1$ factors through $D$.

We extend $\sZ\times_SC_0$ to the cyclic covering $\sZ_C$ of $\BP_C$, obtained by taking the
$d$-th root out of $\sX_C+E$. Remark that $f_D:\sX_D\to D$ is semi-stable. For $\sZ_C\to C$ to
be semi-stable one has to choose $C$ such that the ramification orders over $s_1,\ldots,s_r$ are all
divisible by $2d$. We will assume that this additional condition holds true.

The Higgs bundle $\tau^*(E,\theta)$ has a natural extension to a Higgs bundle
$(\bar E, \bar \theta)$ on $C$ with a Higgs field with logarithmic poles along $C\setminus C_0$.
\begin{claim}\label{8.03}
The pullback of the direct sum decomposition (\ref{eq8.01}) to $C_0$ extends to a decomposition
$$
\big(\bigoplus_{p=0}^{n-1}\bar E^{n-p,p},\bigoplus_{p=0}^{n-1}\bar \theta_{n-p,p}\big)=
\bigoplus_{\iota=1}^\ell\big(\bigoplus_{p=0}^{n-1}\bar F^\iota \,{}^{n-p,p},
\bigoplus_{p=0}^{n-1}\bar \theta^\iota_{n-p,p}\big)
$$
of $(\bar E, \bar \theta)$ over $C$.
\end{claim}
\begin{proof}
By (\ref{sheaves2}) the Hodge bundles of the canonical extension
of $\tau^*(E,\theta)$ are the locally free sheaves
$$
R^pg_{C*}(\Omega^{n-p}_{\BP_C /C}(\log (\sX_C + \Gamma))\otimes \sN^{(i)^{-1}}),
$$
where $\sN^{(i)}=\delta^*(\sO_{\BP(\sE_C)}(1))\otimes \sO_{\BP_C}(-[\frac{iE}{d}])$, and where
$\Gamma$ is the pullback of $C\setminus C_0$.
The decomposition $(R^ng_*\C_{\sZ})_i=\V^1_i\oplus \cdots \oplus \V^\ell_i$
induces a decomposition of the system of Hodge bundles
$$
\bigoplus_{\iota=1}^\ell\big(\bigoplus_{p=0}^{n-1}\bar E_i^\iota \,{}^{n-p,p},
\bigoplus_{p=0}^{n-1}\bar \theta_i^\iota\,{}_{n-p,p}\big),
$$
with $\bar E_i^\iota \,{}^{n-p,p}|_{C_0} = \tau^*E_i^\iota \,{}^{n-p,p}$.

The map $\Phi^\iota_2|_{C_0}$ extends to the cup product map
$$
\bar\Phi_2^\iota:\bar E_i^{\iota \, n-p,p}\otimes \bar p_*(\sN^i) \>>>
R^{p}p_{C*}(\Omega_{\BP_C/C}^{n-p}(\log \sX_C+\Gamma)\otimes \sO_{\BP_C}([\frac{iE}{d}] )).
$$
If $D=C$, since $2i <d$, the divisor $[\frac{iE}{d}]$ is zero.
Let us first show that over $D$ the cup-product
\begin{multline}\label{eq8.04}
\bar\Phi_2:
R^pp_{D*}(\Omega^{n-p}_{\BP_D/D}(\log (\sX_D + \Gamma)\otimes\sN^{-i})
\otimes p_{D*}\sN^i\\
\>>> R^{p}p_{D*}(\Omega_{\BP_D/D}^{n-p}(\log (\sX_D + \Gamma)))
\end{multline}
is surjective. It is sufficient to verify, for $\nu=1, \ldots, i$, the surjectivity of the
product map
\begin{multline*}
R^p p_{D*}(\Omega^{n-p}_{\BP_D /D}(\log (\sX_D + \Gamma))\otimes\sN^{-i+\nu-1})
\otimes p_{D*} \sN\\
\>>> R^{p}p_{D*}(\Omega_{\BP_D/D}^{n-p}(\log (\sX_D + \Gamma))\otimes\sN^{-i+\nu}).
\end{multline*}
For this one can use the pullback under $\delta_D$ of the tautological map
$$
\bigoplus^{n+1}\sO_{\BP(\sE_D)}(-1) \>>> \sO_{\BP(\sE_D)},
$$
and the induced Koszul complex
$$
\sN^{-n-1}\>\alpha_{n+1}>> \bigoplus^{n+1}\sN^{-n} \>\alpha_n>> \cdots \>\alpha_2>>
\bigoplus^{n+1}\sN^{-1} \>\alpha_1>> \sO_{\BP_D}.
$$
Tensorizing with $\Omega^{n-p}_{\BP_D /D}(\log (\sX_D + \Gamma))\otimes \sN^{-i+\nu}$
one has to verify that
\begin{equation}\label{eq8.05}
R^{p'} p_{D*}(\Omega^{n-p}_{\BP_D /D}(\log (\sX_D + \Gamma))\otimes\sN^{-i+\nu-\mu}=0
\end{equation}
for $\mu=1, \ldots, n+1$ and for $p'\geq p+\mu-1$. In fact, this implies
by descending induction on $\mu$, that
$$
R^{p+\mu-1}p_{D*}({\rm Coker}(\alpha_\mu)\otimes \Omega^{n-p}_{\BP_D /D}(\log (\sX_D + \Gamma))
\otimes \sN^{-i+\nu})=0.
$$
Since $0< i-\nu+\mu<d$ the sheaves
\begin{equation}\label{eq8.06}
R^{p'} p_{D*}(\Omega^{n-p}_{\BP_D /D}(\log (\sX_D + \Gamma))\otimes\sN^{(i-\nu+\mu)^{-1}}
\end{equation}
are direct factors of the higher direct
image of $\Omega^{n-p}_{\sZ_D /D}$ (see \cite{EV1}, for example),
hence locally free, and by \ref{1.1} they are trivial. Claim \ref{i1}
allows to assume that $i\leq d-n-2$, hence for $\mu \leq 2$
$$2(i-\nu+\mu)\leq 2(i-\nu+2)\leq 2(i+1)\leq 2d-2n-2 < d.$$
So for $\mu\leq 2$ the sheaves in (\ref{eq8.05}) and (\ref{eq8.06}) coincide.

In general, the sheaves in (\ref{eq8.05}) and (\ref{eq8.06}) differ at most by 
$$H^{p'-1}(E,\Omega^{n-p-1}_E(\log(\sX\cap E))\otimes
\sO_E(E)),$$
which is zero by \ref{1.1}, b), for $p' \geq p+2$.

Remark that the sheaves $R^{p}p_{D*}(\Omega_{\BP_D/D}^{n-p}(\log \sX_D+\Gamma))$
and $p_{D*}\sN^i$ are compatible with base change, whereas
$R^pp_{C*}(\Omega^{n-p}_{\BP_C/C}(\log (\sX_C + \Gamma)\otimes\sN^{(i)^{-1}})$ is
larger than the pullback of the corresponding sheaf from $D$. This is due to the non-semi-stability
of $\sZ_D\to D$ in $p\in C\setminus C_0$.

As the pullback of the direct sum decomposition of the Higgs bundles of the $i$-th eigensheaf
on $D$, we only get a direct sum decomposition $G^1\oplus \cdots\oplus G^\ell$
of a subsheaf of
$$\bigoplus_{p=0}^{n-1}R^{p}p_{C*}(\Omega_{\BP_C/C}^{n-p}
(\log (\sX_C + \Gamma))\sO_{\BP_C}([\frac{iE}{d}] )),$$
containing
$$\bigoplus_{p=0}^{n-1}R^{p}p_{C*}(\Omega_{\BP_C/C}^{n-p}(\log (\sX_C + \Gamma))\sO_{\BP_C}).$$

By \cite{CT}, Section 5, all eigenvalues of the local monodromy
of a Lefschetz pencil for the local system corresponding to the $i$-th eigenspace, are one,
except for one eigenvector $\delta$. Moreover, the local monodromy
is of the form $x\mapsto x\pm h(x)\delta$.

Over $C_0$ we have the decomposition $\tau^*\V^1_i\oplus \cdots \oplus \tau^*\V^\ell_i$. Hence if
the local monodromy is non-trivial, $\delta$ must belong to one of the $\V^\iota_i$, say for
$\iota=1$.
Then for $\iota\neq 1$ the Hodge bundles for $\V^\iota_i$ are contained in the pullback of the
Hodge bundle on the curve $D$, and the image of the direct factor $G^2\oplus \cdots\oplus G^\ell$ in
the direct sum of the cokernels of
$$
R^pp_{C*}(\Omega^{n-p}_{\BP_C/C}(\log (\sX_C + \Gamma)\otimes\sN^{-i})\>\beta_p>>
R^pp_{C*}(\Omega^{n-p}_{\BP_C/C}(\log (\sX_C + \Gamma)\otimes\sN^{(i)^{-1}})
$$
is zero. This implies that
$$
\bigoplus_{p=0}^{n-1}R^pp_{C*}(\Omega^{n-p}_{\BP_C/C}(\log (\sX_C + \Gamma)\otimes\sN^{-i})=
{G'}^1\oplus G^2 \oplus \cdots \oplus G^\ell,
$$
for the subsheaf ${G'}^1=G^1\cap \bigoplus_{p=0}^{n-1}{\rm ker}(\beta_p)$ of $G^1$.
\end{proof}
By \cite{Sim2} $(\bar E,\bar \theta)$ is a poly-stable Higgs bundle of degree zero,
hence \ref{8.03} is a decomposition in sub Higgs bundles, necessarily of degree zero.
Applying \cite{Sim2} again, one obtains a decomposition
\begin{equation}\label{eq8.02}
\tau^* (R^{n-1}f_*(\C_\sX)_{\rm prim})|_{C_0}=\bigoplus_{\iota=1}^\ell \W^\iota_{C_0}
\end{equation}
of $\C$-local systems over $C_0$. Recall that $\sX\times_SC_0$ is isomorphic to the pullback
of the family
$$\Psi:\tilde\BP^n_0=\Psi^{-1}(\BP^1-\{s_1,\ldots,s_r\})\to \BP^1_0=\BP^1-\{s_1,\ldots,s_r\})$$

By Deligne \cite{De3}, I, 1.5, (see also \cite{Voi}, 15.27 and
15.28) the underlying representation of the polarized
variation of Hodge structures $(R^{n-1}\Psi_*\C_{\tilde\BP^n_0})_{\rm
prim}$
$$ \rho:\pi_1(\BP^1_0)\>>> {\rm Sp}(H^{n-1}(F,\Q)_{\rm prim})$$
into the symplectic group with respect to the polarization $<,>$
is irreducible over $\C$. In fact, as explained in \cite{De3}, II,
Section 4.4, one also has the stronger result, that the image of
the monodromy representation is Zariski dense  in ${\rm
Sp}(H^{n-1}(F,\Q)_{\rm prim})$, which is almost simple over $\C.$
Since the Zariski density property of a representation of
$\pi_1(X)$ of an algebraic manifold $X$ in an almost simple
algebraic group remains true under the base change,
the pull back of $\rho$ to a finite covering $C$ of $\BP^1$ is
also irreducible over $\C$.

Hence \ref{8.03} implies that $\ell=1$, contradicting our choice of $i$.
\end{proof}

\begin{example}\label{8.2}
Let  $f:\sX\to S$ be a universal family of hypersurfaces in $\BP^{n}$ of degree $d$.
Let $(g:\sZ\to S)\in \sM_{d,n+1}(S)$ the family of hypersurfaces in $\BP^{n+1}$
obtained by taking the $d$-th root out of $\sX$. By Lemma \ref{2.4}, iii),
$\varsigma(g)=n$ for $d\geq 2(n+1)$. So by \cite{VZ4}, Section 9, the family
$g$ is rigid.
\end{example}
In \ref{8.2} for $d < 2(n+1)$ one has $\varsigma(g)=n-1$, and the family
$$(g:\sZ\to S)\in \sM_{d,n}(S)$$ would be a candidate for a non-rigid family.
Theorem \ref{8.3} says that this is not the case.

\begin{proof}[Proof of Theorem \ref{8.3}]
Let $g:\sZ\to S$ be a universal family for $\sM^{(1)}_{d,n+1}$, obtained
by taking the $d$-th root out of $\sX\to S$. It remains to show, that the family $g:\sZ\to S$
is rigid in $\sM_{d,n+1}$ under the assumption that $n+1 < d < 2(n+1)$.

Assume there exists a morphism $h:\sZ'\to S\times {S'}$, with $S\times {S'} \to M_{d,n+1}$
generically finite and some point $s'_0\in {S'}$ such that the restriction of $h$ to
$h^{-1}(S\times \{s'_0\})$ is isomorphic to $g:\sZ \to S$.
By Proposition \ref{3.3} one has a decomposition
$$ R^nh_*\C_{\sZ'}=\bigoplus_j p_1^*\V^j_S\otimes p_2^*\V_{S'}^j,$$
where $p_i$ denote the projections and where $\V^j_S$ and $\V^j_{S'}$ are irreducible
polarized complex variations of Hodge structures on $S$ and ${S'}$, respectively.
Let us write $E_{S,j}^{q,p}$ and $E_{{S'},j}^{q,p}$ for the corresponding Hodge bundles, and
$F^{n-p,p}$ for the Hodge bundles of $R^nh_*\C_{\sZ}$.
The local Torelli theorem implies that the morphism
\begin{multline*}
T_{S\times {S'}} \>>> {F^{n,0}}^\vee \otimes F^{n-1,1}\\
=\big( \bigoplus_j E_{S,j}^{p,0}\otimes E_{{S'},j}^{n-p,0}\big)^\vee \oplus
\bigoplus_j \big( E_{S,j}^{p-1,1}\otimes E_{{S'},j}^{n-p,0}\oplus
E_{S,j}^{p,0}\otimes E_{{S'},j}^{n-p-1,1} \big)
\end{multline*}
is injective. Then for at least one $j_0$ the restriction
of the Higgs field
$$
(E_{S,j_0}^{p,0}\otimes E_{{S'},j_0}^{n-p,0})\otimes p_2^*T_{S'}  \>>>
E_{S,j_0}^{p,0}\otimes E_{{S'},j_0}^{n-p-1,1}
$$
must be non trivial. This implies that $r={\rm rank}\V_{S'}^{j_0}\geq 2,$
and that $E_{S,j_0}^{p,0}$, $E_{{S'},j_0}^{n-p,0}$ and $E_{{S'},j_0}^{n-p-1,1}$
are all non-zero.
The local system $\V_S^{j_0}$ is irreducible over $\C,$ hence
$\V_S^{j_0}$ is isomorphic to one of the
eigenspaces $(R^ng_{*}\C_{\sZ/S})_{j_0}.$ In particular, \ref{8.0}, ii), implies that $j_0 <d-n-1$.

Restricting everything to $S\times \{s'_0\}\simeq S$ one obtains
$$
\V^{j_0}_S \otimes  \C^r \subset R^nh_*\C_{\sZ'}|_{S\times
\{s'_0\}}=R^ng_*\C_\sZ.
$$
Here $\C^r$ is the trivial local system. As a constant variation of Hodge structures,
it contains non-trivial parts in bidegrees $(n-p,0)$ and $(n-p-1,1)$.
One obtains $r\geq 2$ irreducible factors $\V_1,\ldots,\V_r$ of $R^ng_*\C_{\sZ},$ and up to a shift
in the bidegrees, they are all isomorphic to $(R^ng_{*}\C_{\sZ})_{j_0}.$

By \ref{8.0}, ii), the length of the Griffiths-Yukawa coupling of the latter is $n-1$, and
$\V_\iota$ can only be isomorphic to one of the $(R^ng_{*}\C_{\sZ})_{i(\iota)},$ with
$$i(\iota)\in \{1,\ldots , d-n-1\} \mbox{ \ \ or \ \ }
i(\iota)\in \{n+1,\ldots, d-1\}.$$
In fact, due to the bidegrees of the $\C^r$ part, one can find some $\iota$ and $\iota'$
with
$$i(\iota)\in \{1,\ldots , d-n-1\} \mbox{ \ \ and \ \ }
i(\iota')\in \{n+1,\ldots, d-1\}.$$
However, by \ref{1.6b}, the rank of the first non-vanishing Hodge bundles of the
variations of Hodge structures $(R^ng_{*}\C_{\sZ})_{i(\iota)}$ determines $i(\iota)$,
and $i(\iota)$ must be equal to $j_0$, obviously a contradiction.
\end{proof}

\section{Iterated cyclic coverings}\label{twofold}
Starting with a family $(f: \sX \to S)\in\sM_{d,n}(S)$, we constructed in
\ref{2.1} and \ref{2.2} new families
$$(g_p:\sZ_p \to S) \in \sM_{d,n+p}(S).
$$
For simplicity we will write again $g=g_1$ and $\sZ=\sZ_1$.
By construction, $\sZ_2$ is obtained by taking the $d$-th root out of
$\sZ \subset\BP(\sO_S\oplus \sE)$. By \ref{2.1}
blowing up a section of $\BP(\sO_S\oplus \sE)$, disjoint to $\sZ$,
we obtain an embedding of $\sZ$ to a $\BP^1$ bundle $\BP_1$
and a blowing up $\sY$ of $\sZ_2$ with center a $d$-fold \'etale multisection
of $\sZ_2 \to S$ is obtained by taking the $d$-th root out of $\sZ \subset \BP_1$. The
contraction
$$
\sY \>\psi >> \sZ_2 \subset \BP(\sO_S\oplus \sO_S\oplus \sE)
$$
is again given by the pullback of $\sO_{\BP_1}(1)$ to $\sY$.

Since $\BP_1$ is a $\BP^1$-bundle over $\BP(\sE)$,
there is a morphism $\sY \to \BP(\sE)$ whose fibres are curves.
We will show in this section, that this family of curves is isotrivial
with the Fermat curve $\Sigma_d$ of degree $d$ as general fibre.

Recall that $\Sigma_d$ is the projective curve
with equation $x^d+y^d+z^d=0$, or equivalently, the cyclic cover of $\BP^1$
obtained by taking the $d$-th root out of the divisor of the $d$-th roots of
$1$. We will use different coordinates and represent $\Sigma_d$ as
covering
$$
\Sigma_d \>\alpha'>> \BP^1 \> \beta' >> \BP^1,
$$
where $\beta'$ is the Kummer covering given by $\sqrt[d]{\frac{1}{0}}$
and $\alpha'$ the one given by $\sqrt[d]{\frac{\infty}{{\beta'}^*0}}$.
Writing $G_{\Sigma_d}={\rm Gal}(\Sigma_d/\BP^1)$ we fix a primitive $d$-th root of unit
$\xi_d$ and a generator $\eta_{\Sigma_d}$ of $G_{\Sigma_d}$
which acts on $\sqrt[d]{\frac{\infty}{{\beta'}^*0}}$ by multiplication with
$\xi_d$.

Let $V$ be the normalization of the fibre product $\sZ\times_{\BP(\sE)} \sY$
or, in different terms, the cyclic covering obtained by taking the $d$-th root out
of the pullback of the divisor $\sX\subset \BP(\sE)$ to $\BP_1$.
Writing $H_\bullet=c_1(\sO_\bullet(1))$ in the sequel, this covering is
given by $\sqrt[d]{\frac{\sX}{d\cdot H_{\BP(\sE)}}}.$
Again we choose a generator $\eta_\sZ$ of the Galois group
$G_\sZ={\rm Gal}(\sZ/\BP(\sE))$, acting on $\sqrt[d]{\frac{\sX}{d\cdot H_{\BP(\sE)}}}$
by multiplication with $\xi_d$.

\begin{construction}\label{4.1}
One has a commutative diagram of morphisms between normal varieties
$$
\begin{CD}
V' \>\gamma'>> \sY' \>\tau' >> \BP'_1 \>\mu' >> \BP(\sE)\times \BP^1\\
\A A A \A \delta A A \A \delta_1 A A \A \delta_d A A\\
\hat V \>\hat \gamma >> \hat \sY \>\hat \tau >> \hat \BP_1 \>\hat \mu >> \hat \BP_d\\
\V V V \V \rho V V \V \rho_1 V V \V \rho_d V V \\
V \>\gamma >\sqrt[d]{\frac{\pi^*\sX}{d\cdot \pi^*H_{\BP(\sE)}}}>
\sY \> \tau >\sqrt[d]{\frac{\mu^*E_\sigma}{d\cdot H_{\BP_1}}} > \BP_1 \>\mu >
\sqrt[d]{\frac{E_\infty + d \cdot \pi_d^*H_{\BP(\sE)}}{E_0}} > \BP_d\\
\V V V \V \pi V V \V \pi_1 V V \V\pi_d V V \\
\sZ \>\pi_1>> \BP(\sE) \> = >> \BP(\sE) \> = >> \BP(\sE)
\end{CD}
$$
with:
\begin{enumerate}
\item[(a)] $\delta$, $\delta_1$, $\delta_d$, $\rho$, $\rho_1$ and $\rho_d$
are birational.
\item[(b)] $\pi$ is a family of curves, $\pi_1$ and $\pi_d$
are $\BP^1$ bundles.
\item[(c)] All the vertical arrows are finite coverings
of degree $d$. Here the symbol $\sqrt[d]{f}$ under an arrow indicates, that
the corresponding finite morphism is the Kummer covering given by the $d$-root of $f$.
\end{enumerate}
\end{construction}
\begin{proof}
The morphisms $\tau:\sY \to \BP_1$ and $\BP_1\to \BP_d$ have been considered already
in Section \ref{cyclic}, and $\gamma:V \to \sY$ is just the Kummer covering defined by
$$\sqrt[d]{\frac{\pi^*\sX}{d\cdot \pi^*H_{\BP(\sE)}}}.$$
To complete the construction, we just have to explain
$\delta_d$ and $\rho_d$. Consider on
$\BP_d=\BP(\sO_{\BP(\sE)}\oplus \sO_{\BP(\sE)}(d))$ the tautological
morphism
$$
\pi_d^*(\sO_{\BP(\sE)}\oplus \sO_{\BP(\sE)}(d))\>>> \sO_{\BP_d}(1).
$$
Recall that $E_0$ is the zero set of the section given by the first
direct factor, whereas $E_\sigma$ is given by ${\rm id}\oplus \pi_d^*\sigma$.
Taking the direct sum of both sections, one obtains
$$
\sO_{\BP_d}\oplus \sO_{\BP_d}\>>> \pi_d^*(\sO_{\BP(\sE)}\oplus
\sO_{\BP(\sE)}(d))\>>> \sO_{\BP_d}(1)
$$
and the image of the composite $\Phi$ is $\sO_{\BP_d}(1)\otimes I$, where $I$ is the
sheaf of ideals of $E_0\cap E_\sigma$. The morphism $\rho_d$ is obtained
by blowing up $I$. Let $\hat D_1$ be the exceptional divisor. Then
$(\rho_d\circ \pi_d)^{-1}(\sX)$ is of the form $\hat D_1 + \hat D_2$. The morphism $\delta_d$
is just the blowing down of $\hat D_2$. It is given by the surjection
$$
\pi_d^*(\sO_{\BP(\sE)}\oplus \sO_{\BP(\sE)}) \>>> \rho_d^*(\sO_{\BP_d}(1)(-\hat D_1).
$$
Taking normalizations we obtain the diagram in \ref{4.1}, and a), b) and c).
\end{proof}
\begin{construction}\label{4.2} In \ref{4.1} one has moreover:
\begin{enumerate}
\item[(d)] $V'\cong \sZ\times \Sigma_d$.
\item[(e)] There is a diagonal embedding
$$
G={\rm Gal}(\sZ\times \Sigma_d/\sY')\>>> G'=G_\sZ \times G_{\Sigma_d},
$$
i.e. an embedding whose image is generated by
$(\eta_\sZ,\eta_{\Sigma_d})$. \item[(f)] Taking quotients of $\sZ
\times \Sigma_d$ by $G$ and $G'$ one obtains morphisms
$$\sZ\times \Sigma_d \>\gamma'>>\sY' \> \alpha' >> \BP(\sE)\times \BP^1.
$$
\end{enumerate}
\end{construction}
\begin{proof}
Let us write $\hat E_\bullet$ for the proper transform of $E_\bullet$
under $\rho_d$. Then $\hat\mu$ is the Kummer covering given by
$$
\sqrt[d]{\frac{\hat E_\infty + d \cdot \hat \pi_d^*H_{\BP(\sE)}}{\hat E_0 + \hat D_1}},
$$
where $\hat\pi_d=\rho_d\circ \pi_d$. The morphism $\mu'$ is given by
$$
\sqrt[d]{\frac{\delta_d \hat E_\infty + d \cdot \delta_d\hat \pi_d^*
H_{\BP(\sE)}}{\delta_d\hat E_0 + \delta_d \hat D_1}}=
\sqrt[d]{\frac{\BP(\sE)\times\{\infty \}+ d \cdot H_{\BP(\sE)}\times \BP^1}
{\BP(\sE)\times\{0\} + \sX\times \BP^1}}.
$$
Since $E_0+\pi_d^*(\sX)$ is a normal crossing
divisor, the section $\hat E_\sigma$ does not meet $\hat E_0$ nor $\hat D_2$,
and $\delta_d(\hat E_\sigma)$ neither meets
$$
\delta_d(\hat E_0)=\BP(\sE)\times\{0\}\mbox{ \ nor \ }\delta_d(\hat E_\infty)=
\BP(\sE)\times\{\infty \}.
$$
This allows to choose coordinates in $\BP^1$ with
$\delta_d(\hat E_\sigma)=\BP(\sE)\times\{1 \}$.

Remark that we can choose $H_{\BP_1}$ to be $(\mu^*E_0)_{\rm red}$, hence
$\hat \tau$ is the Kummer covering given by
$$
\sqrt[d]{\frac{\rho_1^*\mu^*E_\sigma}{\rho_1^*\mu^*E_0}}=
\sqrt[d]{\frac{\hat\mu^*\hat E_\sigma}{\hat \mu^* \hat E_0}},
$$
hence $\tau'$ is given by
$$
\sqrt[d]{\frac{{\mu'}^*\BP(\sE)\times\{1 \}}{{\mu'}^*\BP(\sE)\times\{0 \}}}.
$$
Since the last function is the $d$-th root out of the pullback of a function on
$\BP(\sE)\times \BP^1$, one may reverse the order and obtains
\begin{equation}\label{eq4.1}
\sY' \> \alpha' > \sqrt[d]{\frac{\beta^* (\BP(\sE)\times\{\infty \}+ d
\cdot H_{\BP(\sE)}\times \BP^1)}
{\beta^*(\BP(\sE)\times\{0\} + \sX\times \BP^1)}}>
\BP(\sE)\times \BP^1 \> \beta >\sqrt[d]{\frac{\BP(\sE)\times\{1 \}}
{\BP(\sE)\times\{0 \}}}
> \BP(\sE)\times \BP^1.
\end{equation}
Then the composite, considered in f),
\begin{equation}\label{eq4.2}
V' \> \gamma' >\sqrt[d]{\frac{{\alpha'}^*\beta^* (\sX\times \BP^1)
}{{\alpha'}^*\beta^*(d \cdot H_{\BP(\sE)}\times\BP^1)}}> \sY' \> \alpha' >
\sqrt[d]{\frac{\beta^* (\BP(\sE)\times\{\infty \}
+ d \cdot H_{\BP(\sE)}\times \BP^1)}{\beta^*(\BP(\sE)\times\{0\} +
\sX\times \BP^1)}}
> \BP(\sE)\times \BP^1
\end{equation}
is the product of the two cyclic coverings of $\BP(\sE)$ and $\BP^1$,
given by
$$
f_1=\sqrt[d]{\frac{\sX}{d \cdot H_{\BP(\sE)}}}\mbox{ \ \ and \ \ }
f_2=\sqrt[d]{\frac{{\beta'}^*\infty}{{\beta'}^*0}},
$$
respectively, where $\beta':\BP^1 \to \BP^1$ is the $d$-fold cover
ramified at $0$ and $1$. The total space of the first covering is $\sZ$, and
the one of the second is $\Sigma_d$. We obtain d) and an embedding of $G$ into
$G_\sZ\times G_{\Sigma_d}$. For the generators $\eta_\sZ$ and $\eta_{\Sigma_d}$,
acting on $f_1$ and $f_2$ by multiplication with $\xi_d$, the
automorphism $\eta_\sZ\times\eta_{\Sigma_d}$ leaves the function
$$
\sqrt[d]{\frac{\beta^* (\BP(\sE)\times\{\infty \}
+ d \cdot H_{\BP(\sE)}\times \BP^1)}{\beta^*(\BP(\sE)\times\{0\} + \sX\times \BP^1)}}
$$
invariant, hence it generates $G$ as a subgroup of $G_\sZ\times G_{\Sigma_d}$, as
claimed in e).
\end{proof}
One can also reconstruct $\sY$ and the sheaf $\tau^*\sO_{\BP^1}(1)$, defining the
morphism to $\BP(\sO_S\oplus \sO_S\oplus\sE)$ with image $\sZ_2$, starting from the data
in \ref{4.2}. To this aim, write $E'=\BP(\sE)\times {\beta'}^* \infty$,
and let $\zeta:\hat\Pi \to \BP(\sE)\times \BP^1$ be the blowing up of the non singular subscheme
$$
\beta^*(\BP(\sE)\times \{\infty\})\cap \beta^*(\sX \times \BP^1)=
(\BP(\sE)\times E) \cap (\sX\times \BP^1).
$$
Writing $\hat\pi_1:\hat\Pi\to \BP(\sE)$ for the morphism induced
by the projection to $\BP(\sE)$ one has $\hat\pi_1^*(\sX)=\hat
B_1+\hat B_2$, where $\hat B_2$ is the exceptional divisor for
$\zeta$. Let us also write $\hat E$ for the proper transform of
$E'$ under $\zeta$. Blowing down $\hat B_1$ one obtains a morphism
$\eta:\hat\Pi \to \Pi$ to the total space of some projective
bundle $\pi_1:\Pi \to \BP(\sE)$.
\begin{construction}\label{4.3} Using the notations introduced above one has:
\begin{enumerate}
\item[(g)] A finite morphism $\hat\alpha:\hat\sY \to \hat\Pi$ of degree $d$ and
totally ramified over $\hat E + \hat\pi_1^*(\sX)-\hat B_2=\hat E + \hat B_1$, and nowhere else.
In particular $\sY'$ is non-singular.
\item[(h)] A finite morphism $\alpha:\sY \to \Pi$ of degree $d$ with ramification locus
$\eta(\hat E)$.
\item[(i)] Let $\Upsilon=(\hat\alpha^* \hat E)_{\rm red}$. Then the sheaf
$\sO_{\hat \sY}(\Upsilon)\otimes \hat\alpha^*\pi_1^*\sO_{\BP(\sE)}(1)$ defines the
morphism from $\hat\sY$ to $\BP(\sO_S\oplus\sO_S\oplus \sE)$ with image $\sZ_2$.
\end{enumerate}
\end{construction}
\begin{proof}
To verify those properties, let us return to the description of $\sY'$
in (\ref{eq4.1}). Recall that $\hat\BP_d\to \BP(\sE)\times \BP^1$ is given by the
blowing up of
$$
(\BP(\sE)\times \{\infty\})\cap (\sX \times \BP^1),
$$
hence $\hat\Pi$ is finite over $\hat \BP_d$. By (\ref{eq4.2})
$\hat\sY$ is the covering of $\hat \Pi$, given by
$$
\sqrt[d]{\frac{\zeta^*\beta^* (\BP(\sE)\times\{\infty \}+ d \cdot
H_{\BP(\sE)}\times \BP^1)}
{\zeta^*\beta^*(\BP(\sE)\times\{0\} + \sX\times \BP^1)}}=
\sqrt[d]{\frac{\zeta^*(\BP(\sE)\times E' + d\cdot
H_{\BP(\sE)}\times \BP^1)}
{d\cdot \zeta^*(\BP(\sE)\times \{0\}) + \zeta^*(\pi_1^*\sX)}}.
$$
Obviously $\hat B_2$ cancels out, hence it is not in the discriminant
locus, and the latter is non-singular. This implies g), and
since $\hat B_1$ is blown down under $\eta$, part h) as well.

For i) recall that the sheaf $\tau^*\sO_{\BP_1}(1)$ on $\sY$ defines the morphism to $$\BP(\sO_S
\oplus \sO_S\oplus\sE),$$
with image $\sZ_2$. Hence the sheaf $\rho^*\tau^*\sO_{\BP_1}(1)$ is the one we are looking for.

By part g) the divisor $\Upsilon + (\hat\alpha^*\hat B_1)_{\rm red}$ is the reduced
ramification divisor of the covering $\hat \alpha:\hat\sY\to \hat \Pi$. Returning to
the notations in \ref{4.1} one has
$$
\Upsilon=(\hat\tau^*\hat\mu^*(\hat E_\infty))_{\rm red}
\mbox{ \ and \ } (\hat\alpha^*\hat B_1)_{\rm red} = (\hat\tau^*\hat\mu^*\hat D_1)_{\rm red}.
$$
Remark that $\hat\mu^*(\hat E_\infty)_{\rm red}$ is the zero set of a section of
$\rho_1^*(\sO_{\BP_1}(1)\otimes\pi_1^*\sO_{\BP(\sE)}(-1))$, and that both,
$\hat\mu^*(\hat E_\infty)_{\rm red}$ and $\hat\mu^* \hat D_1$ are unramified for
$\hat\tau:\hat\sY\to \hat\BP_1$. Hence
$$
\sO_{\hat\sY}(\Upsilon + (\hat\alpha^*\hat B_1)_{\rm red} )=
\hat\tau^*(\rho_1^*\sO_{\hat\BP_1}(1)\otimes \rho_1^*\pi_1^*\sO_{\BP(\sE)}(-1)\otimes
\sO_{\hat\BP_1}((\mu^* \hat D_1)_{\rm red})),
$$
and one finds
$$
\rho^*\tau^*\sO_{\BP_1}(1)=\sO_{\hat\sY}(\Upsilon)
\otimes \rho^*\pi^*\sO_{\BP(\sE)}(1).
$$
\end{proof}
\section{Deformation of iterated coverings}\label{deformation}
Remark that the description of $\hat\sY$ and of the sheaf defining a morphism
with image $\sZ_2$ in \ref{4.2}, f), and in \ref{4.3}, g) and i), is just using data given
by $\Sigma_d$, by $\sZ$ and by the group action. This allows to replace the Fermat curve
$\Sigma_d$ by any curve $\Sigma$ obtained as a $d$-fold cyclic covering of $\BP^1$, totally
ramified over a reduced divisor $\Gamma$ in $\BP^1$ of degree $d$. The $d$-fold Kummer
covering $\Sigma \to \BP^1$ is given by $\sqrt[d]{f}$ for $f=\frac{\Gamma}{d\cdot 0}$.
We keep the convention, that $\sZ$ is the family of hypersurfaces in $\BP(\sO_S\oplus \sE)$
obtained by taking
the $d$-th root out of a normalized family $\sX \subset \BP(\sE)$. We fix a generator
$\eta_\sZ$ of the Galois group, acting on $\sqrt[d]{\frac{\sX}{d\cdot H_{\BP(\sE)}}}$
by multiplication with $\xi_d$.

One has a $2d$-fold covering $\sZ\times \Sigma \to \BP(\sE)\times \BP^1$ with Galois group
$G_\sZ\times G_\Sigma$. We choose a generator $\eta_\Sigma$, acting on
$\sqrt[d]{f}$ by multiplication with $\xi_d$ and we define $G\subset G_\sZ\times G_\Sigma$
to be the subgroup generated by $(\eta_\sZ,\eta_\Sigma)$.

By abuse of notations we will not add any additional index to the varieties and maps
considered in the last section, keeping however in mind, that the diagram
in \ref{4.1} does not exists in general. Let us collect the surviving properties:
\begin{lemma}\label{5.1} Let $E=\BP(\sE)\times \Gamma$. One has a commutative diagram
$$
\begin{CD}
\sZ\times \Sigma \> \gamma' >> \sY' \> \alpha' >> \BP(\sE)\times \BP^1\\
\noarr \A \delta A A \A A \zeta A \\
\noarr \hat \sY \> \hat \alpha >> \hat \Pi \>\hat \pi >> \BP(\sE)
\end{CD}
$$
where
\begin{enumerate}
\item[(i)] $\gamma'$ is the quotient by $G$ and $\alpha'$ is the quotient by $G'/G$.
\item[(ii)] $\zeta$ is the blowing up of $E\cap(\sX\times \BP^1)$ with exceptional divisor
$\hat B_2$.
\item[(iii)] $\hat\alpha$ is the finite covering, totally ramified over $\hat E + \hat B_1$,
where $\hat E$ and $\hat B_1$ denote the proper transforms of $E$ and $\sX\times \BP^1$, respectively. In particular $\hat \sY$ is non singular.
\item[(iv)] $\delta$ is the blowing up of a non-singular subvariety, consisting of $d$
copies of $\sX$. The exceptional divisor $\Delta_2$ is a $\BP^1$ bundle over
$d$ disjoint copies of $\sX$.
\end{enumerate}
\end{lemma}
\begin{proof} Obviously $\sZ\times \Sigma/G'=\BP(\sE)\times \BP^1$, and one obtains i).
In ii) remark that $E\cap (\sX\times \BP^1)$ consists of $d$ disjoint copies of $\sX$.
For iii) remark that $\alpha':\sY'\to \BP(\sE)\times \BP^1$ is the Kummer covering defined
by
$$
\sqrt[d]{\frac{E + d\cdot H_{\BP(\sE)}\times \BP^1}{d\cdot \BP(\sE)\times 0 +
\sX\times \BP^1}}.
$$
So $\hat B_2$ is not part of the ramification locus. $\Delta_2$ is the preimage of $\hat B_2$
hence a covering of a $\BP^1$-bundle, ramified along two disjoint sections. One obtains iv).
\end{proof}

\begin{corollary}\label{5.2}
If in \ref{5.1} $\Sigma$ is the Fermat curve $\Sigma_d$, then $\hat\alpha:\hat\sY \to
\hat\Pi$ coincides with the morphism considered in \ref{4.3}, g).
\end{corollary}
It remains to reconstruct the non-singular model $\sZ_2$ of $\hat\sY$, and the
invertible sheaf defining the contraction from $\hat\sY \to \sZ_2$.
Taking \ref{4.3} as a model, and using the notation from \ref{5.1}, we define
the divisor $\Upsilon=(\hat\alpha^*\hat E)_{\rm red}$ and the invertible sheaf
$$
\sN=\sO_{\hat\sY}(\Upsilon)\otimes \hat\pi^*\sO_{\BP(\sE)}(1).
$$
We will write $\hat g_2=p\circ\hat\pi\circ\hat\alpha: \hat\sY \to S$,
and $\Delta_i=\hat\alpha^*B_i$
\begin{lemma}\label{5.3} \ \
\begin{enumerate}
\item[(i)] The sheaf $\sN$ is generated by global sections.
\item[(ii)] $\hat g_{2*} \sN= \sE \oplus \sO_S^{\oplus 2}$.
\item[(iii)] The image of the morphism $\phi:\hat\sY \to \BP(\hat g_{2*} \sN)$,
defined by $\hat g_2^*\hat g_{2*} \sN$, is non-singular.
\item[(iv)] $\phi|_{\hat\sY\setminus (\Upsilon+\Delta_1)} $ is an embedding,
$\phi|_{\Delta_1}$ is a $\BP^1$ bundle over $\sX$ and $\phi$ contracts
$\Upsilon$ to a section of $\BP(\hat g_{2*}\sN) \to S$.
\item[(v)] $\phi$ factors as
$$\hat\sY \> \rho >> \sY \> \psi >> \sZ_2,$$
where $\sY\to S$ and $g_2:\sZ_2\to S$ are smooth and where $\rho$ contracts $\Delta_1$
to $\sX\subset \sY$ and where $\psi$ contracts $\rho(\Upsilon)$ to a section
of $g_2$.
\end{enumerate}
\end{lemma}
\begin{proof}
Let us write $\sO_{\hat\Pi}(i,j)$ for the pullback of
$pr_1^*\sO_{\BP(\sE)}(i)\otimes pr_2^*\sO_{\BP^1}(j)$.
By definition
\begin{equation}\label{eq5.1}
\sN=\sO_{\hat\sY}(\Upsilon)\otimes \hat\alpha^*\sO_{\hat\Pi}(1,0),
\end{equation}
hence it is globally generated outside of $\Upsilon$.
On the other hand one has
$$
\sO_{\hat\sY}(\Upsilon-\Delta_1)=
\sO_{\hat\sY}((\Upsilon+\Delta_2)-(\Delta_1+\Delta_2)) = \hat\alpha^*\sO_{\hat\Pi}(1,-1),
$$
and hence
\begin{equation}\label{eq5.2}
\sN=\sO_{\hat\sY}(\Delta_1)\otimes \hat\alpha^*\sO_{\hat\Pi}(0,1).
\end{equation}
Then $\sN$ is also globally generated outside of $\Delta_1$,
and since $\Delta_1\cap \Upsilon=\emptyset$ i) holds true.

The covering $\hat\alpha$ is obtained by taking the $d$-th root out of
$$\hat E + (d-1)\cdot \hat D_1 + d \cdot \hat D_2.$$
The corresponding invertible sheaf is the pullback $\sO_{\hat\Pi}(d(d-1),d)$ of
$$pr_1^*\sO_{\BP(\sE)}(d(d-1))\otimes pr_2^*\sO_{\BP^1}(d).$$
As in Section \ref{cyclic} this implies that
\begin{align*}
\hat\alpha_*\sO_{\hat\sY}&=\bigoplus_{i=0}^{d-1}\sO_{\hat\Pi}(-i(d-1),-i)
(i \cdot \hat D_2 + [\frac{i(d-1)}{d}] \cdot \hat D_1)\\ &=
\sO_{\hat\Pi}\oplus \bigoplus_{i=1}^{d-1}\sO_{\hat\Pi}(i-d,-i)
(\hat D_2 ),
\end{align*}
and
\begin{align}\notag
\hat\alpha_*\sN &=
\hat\pi_1^*\sO_{\BP(\sE)}(1)\oplus \sO_{\hat\Pi}(-1,-d+1)
(\hat E + \hat D_2 )\otimes \hat\pi_1^*\sO_{\BP(\sE)}(1)\\
\label{eq5.3}
&\oplus \bigoplus_{i=1}^{d-2}\sO_{\hat\Pi}(i-d,-i)
(\hat D_2 )\otimes \hat\pi_1^*\sO_{\BP(\sE)}(1)\\ \notag
&=\sO_{\hat\Pi}(1,0) \oplus \sO_{\hat\Pi}(0,1)
\oplus \bigoplus_{i=1}^{d-2}\sO_{\hat\Pi}(1+i-d,-i)
(\hat D_2 ).
\end{align}
In this decomposition, applying $\hat\pi_{*}$ to the right hand side
of this decomposition, one obtains zero. Hence
$$
\hat g_{2*} \sN = \hat\pi_{*}(\sO_{\hat\Pi}(1,0) \oplus \sO_{\hat\Pi}(0,1)).
$$
We obtain part ii), as well as the existence
of the morphism $\phi$.

The sheaf $\hat g_{2*}\sN$ is locally free, and by ``cohomology and base change''
one may assume in iii) that $S$ is a point.
By (\ref{eq5.1}) and (\ref{eq5.2}) outside
of $\Upsilon\cup\Delta_1$ the morphism $\phi$ factor through
an \'etale morphism, the one to $\BP(\sE)\times \BP^1$. On the other hand,
it separates points on each fibre of $\hat\pi_1$, hence of $\hat\alpha$, hence
it is an embedding on the complement of $\Upsilon\cup\Delta_1$.

The morphism $\hat\alpha_* \sN(-\Upsilon) \to  \hat\alpha_* \sN$
induces an isomorphism for all of the direct factors in (\ref{eq5.3}),
except of the second one. Therefore $\sN|_\Upsilon=\sO_\Upsilon$. By construction
$\Upsilon$ is isomorphic to $\BP^{n}$ and since $\Upsilon + \Delta_2$ is a fibre
of a morphism, this also implies that $\sO_\Upsilon(\Upsilon)=\sO_{\BP^n}(-1)$.
So $\phi$ contracts $\Upsilon$ to a point. Since for each point
$p\in \Upsilon$ the restriction of $\phi$ to a transversal curve is an embedding,
the image point of $\Upsilon$ is a non-singular point of $\phi(\hat\sY)$.

As in (\ref{eq5.1}) and in (\ref{eq5.2}) we can reverse the role of $\Upsilon$ and
$\Delta_1$, and repeating the arguments used above one finds in this case
that $\sN|_{\Delta_1}=\hat\alpha^*\sO_{\hat\Pi}(1,0)|_{\Delta_1}$.
So $\phi: \Delta_1 \to \phi(\Delta_1)$
is a $\BP^1$ bundle over $\sX$, with $\sO_{\Delta_1}(\Delta_1)$ fibrewise isomorphic to
$\sO_{\BP^1}(-1)$, and $\phi(\hat\sY)$ is smooth in a neighborhood of
$\phi(\Delta_1)$.
\end{proof}
The moduli scheme $M_{d,1}$ of $d$ distinct points in $\BP^1$ is irreducible and of
dimension $d-3$, for $d \geq 3$. Let $T\to M_{d,1}$ be a generically finite
and surjective morphism from a non-singular projective variety $T$.
Assume that $T\to M_{d,1}$ is induced by a
normalized family $g_0:\sW_0 \to T$ of points in $\BP^1$.
So for some vectorbundle $\sF$ of rank $2$ there is an embedding $\sW_0 \to \BP(\sF)$,
with $\sO_{\BP(\sF)}(d)=\sO_{\BP(\sF)}(\sW_0)$.

By \ref{2.1} we obtain a family $g:\sW \to T$ of non-singular
degree $d$-curves in $\BP(\sO_T\oplus \sF)$. By assumption, one of
the fibres is the Fermat curve $\Sigma_d$, say the one over
$t_0\in T$. Let $\sX \to S$ be a universal family for $\sM_{d,n},$ and
let $(\sZ \to S) \in \sM_{d,n}(S)$ be the family of cyclic
covers, obtained in \ref{2.1}. For each point $t\in T$ \ref{5.1}
and \ref{5.3} gives a family $\sZ_{2t} \to S\times
\{t\}$. Of course, the explicit construction of this
family in \ref{5.1} extends to families of curves $\Sigma$, and
\ref{5.3} extends by base change.

So there exists a family $(\sZ_2 \to S\times T)\in
\sM_{d,n+2}(S\times T)$. Obviously, a general fibre of this family has
only one automorphism of degree $d$, hence the curve $\Sigma$ as well as the
corresponding fibre in $\sZ$ are uniquely determined. One obtains:
\begin{proposition}\label{5.4} For $n\geq 3$ and $d\geq n+1$
consider the sub-moduli stack $\sM_{d,n}^{(2)}$ of $\sM_{d,n},$
and a universal family $g_2:\sZ_2 \to S_2$ for $\sM^{(2)}_{d,n}$.
Then there exists a $(d-3)$-dimensional manifold $T$, and generically
finite morphisms
$$\varphi': T \>>> M^{(1)}_{d,2}\mbox{ \ \ and \ \ }
\varphi:S=S_2\times T \>>> M_{d,n}$$
with:
\begin{enumerate}
\item[(a)] $\varphi$ is induced by a normalized family
$g:\sZ \to S$.
\item[(b)] Let $t_0\subset T$ be a point whose image in $M_{d,1}$
is the moduli point of the Fermat curve. Then the restriction of $g$ to
$S\times \{t_0\}$ coincides with $h_2$.
\item[(c)] If for some $t\in T$ the point $\varphi'(t)$ is the moduli
point of a curve $\Sigma$, then $g^{-1}((s,t))$ is the quotient of
$g_2^{-1}(s)\times \Sigma$ by the action of $\Z/d$.
\end{enumerate}
\end{proposition}

\begin{proof}[Proof of Theorem \ref{5.5}]
By induction we may assume that there is a generically finite
morphism $S_\nu \times T^{\times (r-1)} \to M_{d,n-2}$, which is induced by a family.
Then applying \ref{5.4} one obtains a generically finite morphism
$$
S_\nu \times T^{\times r} \>>> M^{(2)}_{d,n}\times T \>>> M_{d,n},$$
again induced by a family.
\end{proof}
\section{Variation of Hodge structures for iterated coverings}\label{variation}

Keeping the notations from Section \ref{twofold} we will compare the
variations of Hodge structures of the normalized family
$$(g_2:\sZ_2\to S) \in \sM_{d,n+2}(S),$$
with the one of $(f:\sX \to S)\in \sM_{d,n}(S).$
As it will turn out, this can be done using only the properties
stated in \ref{4.2}, and \ref{4.3}, and the results carry over to the slightly
more general situation considered in \ref{5.1} and \ref{5.3}.

Recall that in \ref{4.1} or \ref{5.3}, v), we considered a blowing up $\psi:\sY\to\sZ_2$
with center a section of $g_2$. So obviously on has
\begin{claim}\label{6.1} There exists a constant $\Q$ variation of Hodge structures
$\W'$ with
$$R^{n+1}(g_2\circ \psi)_*\Q_{\sY}=\W'\oplus R^{n+1}g_{2*}\Q_{\sZ_2}.$$
Moreover, $\W'$ is trivial for $n$ even and of rank one, concentrated in bidegree
$(\frac{n+1}{2},\frac{n+1}{2})$ for $n$ odd.
\end{claim}
Let us compare next the variations of Hodge structures for
$\sY'$, $\hat\sY$ and $\sY$. Remark that $\sY'$ fibrewise
only has quotient singularities, hence the Hodge structures of the fibres of
the morphism $g'_2:\sY'\to S$ are pure.
In the sequel $(-j)$ denotes the Tate twist of a variation of Hodge structures,
i.e. the shift of the bigrading by $(j,j)$.
\begin{claim}\label{6.2}
The morphisms $\rho:\hat\sY \to \sY$ and $\delta:\hat\sY \to \sY'$ induce Hodge
isometries
$$ R^{n+1}(g_2\circ\psi\circ\rho)_*\Q_{\hat\sY}\simeq
R^{n+1}(g_2\circ \psi)_*\Q_{\sY}\oplus
R^{n-1}g_{2*}\Q_{\sX}(-1)$$
and
$$
R^{n+1}(g_2\circ\psi\circ\rho)_*\Q_{\hat\sY}\simeq
R^{n+1}g'_{2*}\Q_{\sY'}\oplus \bigoplus^d
R^{n-1}g_{2*}\Q_{\sX}(-1).$$
Moreover,
$$
R^{n+1}(g_2\circ \psi)_*\Q_{\sY}\simeq
R^{n+1}g'_{2*}\Q_{\sY'}\oplus \bigoplus^{d-1}
R^{n-1}g_{2*}\Q_{\sX}(-1).
$$
\end{claim}
\begin{proof}
The first two equalities follow from the explicit description in \ref{5.3}, v), and
\ref{5.1}, iv), of $\rho$ and $\delta$ as blowing up with centers isomorphic to $\sX$ and to
$d$ copies of $\sX$, respectively. For the last one remark, that the exceptional
divisors $\Delta_1$ of $\rho$ and $\Delta_2$ of $\delta$ meet transversally
in $d$ sections of the $\BP^1$-bundle $\Delta_1 \to \sX$.
\end{proof}
Recall next, that $\sY'=\sZ\times \Sigma/G$, where $G\cong \Z/d$
is diagonally embedded in $G_\sZ\times G_{\Sigma}$.
Write $H^b(\Sigma,\Q_{\Sigma})_{i}$ for the direct factor of the constant
variation of $\Q(\xi_d)$-Hodge structures $H^b(\Sigma,\Q(\xi_d)_{\Sigma})$,
where the generator $\eta_{\Sigma}$ of $G_{\Sigma}$
acts by multiplication with $\xi_d^i$, and correspondingly
$(R^{a}g_{*}\Q_{\sZ})_i$ for the sub-variation of $\Q(\xi_d)$-Hodge structure
of $R^{a}g_{*}\Q(\xi_d)_{\sZ}$ where $\eta_\sZ$ acts in the same way.
Remark that
$$(R^{a}g_*\Q_{\sZ})_i\otimes H^b(\Sigma,\Q_{\Sigma})_{d-i}
$$
has a $\Q$-structure. Obviously one obtains:
\begin{claim}\label{6.3}
\begin{multline*}
R^{n+1}g'_{2*}\Q_{\sY'}=\big(\bigoplus_{a+b=n+1} R^a g_* \Q_\sZ \otimes
H^b(\Sigma,\Q_{\Sigma})\big)^{G}=\\
\bigoplus_{a+b=n+1}\bigoplus_{i=0}^{d-1} (R^a g_* \Q_\sZ)_i \otimes
H^b(\Sigma,\Q_{\Sigma})_{d-i}.
\end{multline*}
\end{claim}
Of course in the decomposition in \ref{6.3} one only has to consider
$b=0,1$ and $2$. Since
$$H^0(\Sigma,\Q_{\Sigma})_{d-i}=H^2(\Sigma,\Q_{\Sigma})_{d-i}=0$$
for $i\neq 0$, and since
$$(R^a g_* \Q_\sZ)_0=R^a p_* \Q_{\BP(\sO\oplus\sE)}=0$$
for $a$ odd, the direct factors for $i=0$ in the decomposition in Claim \ref{6.3}
vanish for $n$ even. For $n$ odd, they are constant variation of Hodge structures,
concentrated in bidegree $(\frac{n+1}{2},\frac{n+1}{2})$. One obtains:
\begin{proposition}\label{6.4}
There exist constant $\Q$ variations of Hodge structures $\W'$ and $\W$ and a
Hodge isometry
\begin{multline*}
R^{n+1}g_{2*}\Q_{\sZ_2}\oplus \W' \simeq
\bigoplus_{1}^{d-1} (R^n g_* \Q_\sZ)_i \otimes
H^1(\Sigma,\Q_{\Sigma})_{d-i}
\oplus \bigoplus^{d-1} R^{n-1}f_*\Q_{\sX}(-1)
\oplus \W.
\end{multline*}
If $n$ is even, $\W'=\W=0$. If $n$ is odd, $\W$ and $\W'$ are concentrated in
bidegree $(\frac{n+1}{2},\frac{n+1}{2})$.
\end{proposition}
The Proposition \ref{6.4}, applied to $\Sigma=\Sigma_d$, allows to complement Corollary
\ref{2.5}.

\begin{proof}[Proof of Theorem \ref{6.5}]
Part b) has been shown in \ref{2.5}. For a)
start with a family $f:\sX \to S$, coming from a dominant and generically finite
morphism $S\to M_{d,n-\ell}$. Replacing $S$ by some covering, we may assume that the family
is normalized. Since
$$\ell \geq n-[\frac{d}{2}]+1, \mbox{ \ \ hence \ \ }d\geq 2(n-\ell+1),$$
\ref{2.4}, iii), implies that for the $d$-fold covering family
$\varsigma(g_1)=n-\ell$. Using the decomposition in \ref{6.4} one finds
$\varsigma(g_2)=n-\ell$.
Since $n-\ell$ is strictly smaller than the dimension of the
fibres, one can apply \ref{2.4}, iv), and the length of the Griffiths-Yukawa coupling
remains $n-\ell$ for all families, obtained by further $d$-fold coverings.
\end{proof}

\section{Quintic threefolds with complex multiplication}\label{quintic}
Recall the definition of the special
Mumford-Tate group (called Hodge group in \cite{Mum1} and \cite{Mum2}),
Let $V$ be a $\Q$-vector space with a Hodge structure of weight $k$.
By \cite{De3} the Hodge structure on $V$ is given by a homomorphism of real algebraic groups
$$h: S\>>> {\rm Gl}(V)_{\R},$$
where $S$ is the real algebraic group obtained from $\C^*$ by restriction of scalars from $\C$
to $\R$. Let $U^1$ denote the subgroup $U^1=\{z\in \C^*;z\bar z =1\}$.
The special Mumford-Tate group $\Hg(V)=\Hg(V,h)$ is the smallest algebraic subgroup of
${\Gl}(V)$ defined over $\Q$, with
$$
h(U^1) \subset \Hg(V,h)_\R=\Hg(V,h)\otimes_\Q\R.
$$
It is also the largest $\Q$ algebraic subgroup of $\Gl(V)$,
which leaves all Hodge tensors invariant, i.e. all elements
$$
\eta\in \big[ \big(\bigotimes^{m}V\big) \otimes
\big(\bigotimes^{m'}V^\vee\big) \big]^{\frac{k(m-m')}{2},
\frac{k(m-m')}{2}}.
$$
For a manifold $F$ consider the Hodge structure $V=H^{\dim(F)}(F,\Q)_{\rm prim}$. The special
Mumford-Tate group of $F$ is $\Hg(F)=\Hg(V)$, necessarily a reductive group.

One says that $F$ (or a $\Q$ Hodge structure) has complex multiplication,
if $\Hg(H^{\dim F}(F,\Q)_{\rm prim},h)$ (or $\Hg(V)$) is a commutative group. One also says that
$F$ (or $V$) is of CM type in this case.

Note that for a Calabi-Yau $3$-fold $F$
$$H^3(F,\Q)_{\rm prim}=H^3(F,\Q),$$
as $H^1(F,\Q)=0$, hence that $\Hg(F)=\Hg(H^3 (F,\Q)).$ If $F$ is a surface, again
$\Hg(F)=\Hg(H^2 (F,\Q))$ since ${\rm NS}(F)\otimes \Q$ is invariant under $\Hg(F)$.

\begin{lemma}\label{7.1} \
\begin{enumerate}
\item[(a)] If $V$ and $W$ are two $\Q$-Hodge structures of weight $k$, then
$$\Hg(V\oplus W)\subset
\Hg(V)\times \Hg(W) \subset \Gl(V)\times \Gl(W) \subset \Gl(V\oplus W),$$
and the projections
$$\Hg(V\oplus W)\>>> \Hg(V),\mbox{ \ \ and \ \ } \Hg(V\oplus W)\>>> \Hg(W)$$
are surjective.
\item[(b)] The special Mumford-Tate group does not change under Tate twists,
i.e. $\Hg(V(-1))=\Hg(V)$.
\item[(c)] The special Mumford-Tate group of a Hodge structures concentrated
in bidegree $(p,p)$ is trivial.
\item[(d)] If $V$ and $W$ are two $\Q$ polarized Hodge structures, then
$V\otimes W$ has complex multiplication, if and only if both, $V$ and $W$
have complex multiplication.
\end{enumerate}
\end{lemma}
\begin{proof}
Assume in a) that the Hodge structures on $V$ and $W$ are given by
$$
h_1:S \>>> \Gl (V)_\R \mbox{ \ \ and \ \ } h_2:S \>>> \Gl (W)_\R,
$$
respectively. Obviously
$$\Hg(V, h_1)_\R \times \Hg(W, h_2)_\R \subset \Gl(V)_\R\times \Gl(W)_\R \subset \Gl(V\oplus W)_\R$$
is defined over $\Q$ and it contains the image of the homomorphism
$$
(h_1,h_2): U^1 \>>> {\rm GL}(V)_\R\times {\rm GL}(W)_\R\subset \Gl(V\oplus W)_\R.$$
By definition this implies that
$$\Hg(V\oplus W, h_1\oplus h_2)\subset\Hg(V, h_1)\times \Hg(W, h_2).$$
The image of the projection
$$\Hg(V\oplus W, h_1\oplus h_2)\subset\Hg(V, h_1)\times \Hg(W, h_2)\>>> \Hg(V, h_1)$$
is an $\Q-$algebraic subgroup $G\subset \Hg(V, h_1).$ Since $G_\R$ contains the image of
$h_1:U^1\to {\rm GL}(V)_\R$ one finds $G=\Hg(V,h_1).$

The parts (b) and (c) are obvious, and d) is Prop. 1.2 in \cite{Bor}.
\end{proof}
Let $f: \sX\to S$ be a universal family of five points in $\BP^1,$
and let $g_1: \sZ_1\to S$ be the family of the $5$-th cyclic
covers of $\BP^1$ ramified on $\sX.$ Note that this family is one
of the example in \cite{DM}. The Galois group $G=\Z/5$ acts
fibrewise by automorphisms on the family $g_1:\sZ_1\to S$. We
consider the induced family $g_1: {\rm Jac}(\sZ_1/S)\to S$ of
Jacobians.  Let  $\xi=e^{\frac{2\sqrt{-1}\pi}{5}}.$ Then $\Z(\xi)$
acts as a sub-ring of the endomorphism ring of  ${\rm
Jac}(\sZ_1/S)\to S$ via the action of $\xi$ on $g_1: \sZ_1\to S.$ The
intersection form $<,>$ on the $\Q$-variation of Hodge structures
$$R^1g_{1,*}\Q_{\sZ_1}=R^1g_{1,*}\Q_{{\rm Jac}(\sZ_1)}$$
is defined by taking cup product of 1-formes along the fibres of $g_1: \sZ_1\to S$.

\begin{claim}\label{7.2}
For $l\in \Z(\xi)$ and for all $x,\,y\in R^1g_{1,*}\Q_{\sZ_1}|_{s_0}=H^1(g^{-1}_1(s_0),\Q)$
one has $<lx,y>=<x, \bar ly>$.
\end{claim}
\begin{proof}
In fact, the corresponding property holds true for all cyclic coverings of degree $n$.
Let $\sigma$ be a generator of $G=\Z/n.$ Then
$$<\sigma x,\sigma y>=<x,y>,\quad \forall x,\, y\in H^1(g^{-1}_1(s_0),\Q).$$
Let
$$ H^1(g^{-1}_1(s_0),\Q)=\bigoplus_{i=1}^{n-1}V_i,$$
be the decomposition in eigenspaces, i.e. $\sigma(v)=\xi^iv$ for all $v\in V_i.$
Then $<V_i,V_j>=0$ for all $i, j$ with $i+j\not=n.$
On the other hand, for $x\in V_i$ and $y\in V_{n-i}$, the equality
$\bar\sigma=\sigma^{-1},$ implies that
\begin{multline*}
 <\sigma x, y>=\xi^i<x,y>=<x,\xi^iy>=<x,(\xi^{n-i})^{-1}y>\\
 =<x,\sigma^{-1}y>=<x,\bar\sigma y>.
 \end{multline*}
\end{proof}
We are therefore in the situation of \cite{De2}, 4.9.
The following construction is similar to \cite{DeJN} In fact, the second family in \cite{DeJN}
lies in the degeneration of $g_1:\sZ_1\to S$.

As in \cite{DeJN}, one starts with the group $G$, constructed by
Deligne. For a compact open subgroup $K\subset G({\bf A}_f)$ the
quotient ${_K}M_\C(G,h_0)$ is the moduli space of isomorphism
classes of principally polarized Abelian varieties of dimension
$6,$ together with the given $\Z(\xi)$-action satisfying the
property in Claim \ref{7.2} and a level 1 structure (\cite{De2},
4,12 and \cite{DeJN}, Section 2). For a suitable choice of $K$,
the family of Jacobians induces a generically finite  morphism
$$\phi: S\>>> {_K}M_\C(G,h_0).$$

\begin{claim}\label{7.3}
The map $\phi: S\to {_K}M_\C(G,h_0)$ is dominant.
\end{claim}
\begin{proof}
Since $\dim S=2$, and since $\phi$ is generically finite, we only need to show that
$$\dim ({_K}M_\C(G,h_0))=2.$$
After base change we may assume that there exists a universal
family
$$ \pi: \mathcal A\>>> {_K}M_\C(G,h_0),$$
together with an $\Z(\xi)$-action on the fibres.
This leads the eigenspace decomposition of $R^1\pi_*(\Q_{\mathcal A})$ as polarized complex
variation of Hodge structures
$$R^1\pi_*(\Q_{\mathcal A})\otimes\Q(\xi)=\V(\xi)\oplus\V(\xi^2)\oplus\V(\xi^3)\oplus
\V(\xi^4).$$
Since $<lx, y>=<x,\bar ly>,$ the intersection form $<,>$ induces a perfect duality between
$\V(\xi^i)$ and $\V(\xi^{5-i}).$\\

Next we determine the ranks of the Hodge bundles in the corresponding decomposition.
Note that the pull back of $R^1\pi_*(\Q_{\mathcal A})$ together with the $\Z(\xi)$-action is
just $R^1g_{1,*}(\Q_{{\rm Jac}(\sZ_1)})$ together with the $\Z(\xi)$-action. We only need to
determine the ranks of the Hodge bundles in the decomposition
$$ R^1g_{1,*}(\Q_{{\rm Jac}(\sZ_1)})\otimes\Q(\xi)=\W(\xi)\oplus\W(\xi^2)\oplus\W(\xi^3)\oplus
\W(\xi^4).$$
Writing $h^{p,q}(\xi^i)$ for the rang of the $(p,q)$ Hodge bundle of $\W(\xi^i)$ one has
$$h^0(\Omega^1_{\BP^1}(5-i))=h^{1,0}(\xi^i),\mbox{ \ \ and \ \ } h^1(\BP^1,\sO(-i))=h^{0,1}(\xi^i).$$
Then one finds
\begin{gather*}
(h^{1,0}(\xi), h^{0,1}(\xi))=(3,0),\ \ \ \ \ \ (h^{1,0}(\xi^2), h^{0,1}(\xi^2))=(2,1),\\
(h^{1,0}(\xi^3), h^{0,1}(\xi^3))=(1,2), \mbox{ \ \ and \ \ }(h^{1,0}(\xi^4), h^{0,1}(\xi^4))=(0,3).
\end{gather*}
In particular, $\V(\xi)$ and $\V(\xi^4)$ are unitary local subsystems. The perfect duality between
$\V(\xi^2)$ and $\V(\xi^3)$  implies that the corresponding  the Higgs bundles
$$ E^{1,0}(\xi^2)\>>> E^{0,1}(\xi^2)\otimes \Omega^1_{{_K}M_\C(G,h_0)},\quad E^{1,0}(\xi^3)
\>>> E^{0,1}(\xi^3)\otimes \Omega^1_{{_K}M_\C(G,h_0)}$$
are dual to each other.  This gives a precise description of the rang of the differential
map
$$ d: T_{{_K}M_\C(G,h_0)}\>>>  S^2E^{0,1}\subset {E^{1,0}}^{\otimes 2}$$
of the natural inclusion of ${_K}M_\C(G,h_0)$ into the moduli
space of the polarized Abelian varieties in terms of the above
eigenspace decomposition.  Since $\V(\xi),\,\V(\xi^4)$ are
unitary, the differential map $d$ factors over
$$ d: T_{{_K}M_\C(G,h_0)}\>>> (E^{1,0}(\xi^2)\oplus E^{1,0}(\xi^3))^{\vee}
\otimes(E^{0,1}(\xi^2)\oplus E^{1,0}(\xi^3)).$$
Since the Higgs field preserves the eigenspace decomposition
$$ (E^{1,0}(\xi^2)\oplus E^{0,1}(\xi^2))\oplus (E^{1,0}(\xi^3)\oplus E^{0,1}(\xi^3)),$$
and $d$ factors further  through the diagonal map
\begin{multline*}
d: T_{{_K}M_\C(G,h_0)}\>>> {E^{1,0}(\xi^2)}^{\vee}\otimes E^{0,1}(\xi^2)\oplus
{E^{1,0}(\xi^3)}^{\vee}\otimes E^{0,1}(\xi^3)\\
\simeq ({E^{1,0}(\xi^2)}^{\vee}\otimes E^{0,1}(\xi^2))^{\oplus 2}.
\end{multline*}
The generical injectivity of $d$ implies that the Kodaira-Spencer map on the each copy
$$\theta_{1,0}: T_{ {_K}M_\C(G,h_0)}\>>> {E^{1,0}(\xi^2)}^{\vee}\otimes E^{0,1}(\xi^2)$$
also is injective. Hence,
$$2\leq \dim({_K}M_\C(G,h_0))\leq {\rm rank}({E^{1,0}(\xi^2)}^{\vee}\otimes E^{0,1}(\xi^2))=2.$$
\end{proof}
\begin{corollary}\label{7.4}
Let $S'$ be the set of $s\in S$ such that $g_1^{-1}(s)$ has complex multiplication. Then
$S'$ is dense in $S$.
\end{corollary}
\begin{proof}
By (\cite{Mum2}, Section 2, \cite{De2}, 5.1 and 5.2) the set of CM points in ${_K}M_\C(G,h_0)$
is dense.
\end{proof}
Consider now the second and third iterated $5$-fold coverings
$g_2:\sZ_2 \to S$ and $g_3:\sZ_3\to S$. Replacing $S$ by an \'etale covering,
we may assume that $F:\sX \to S$ consists of $5$ disjoint sections, hence
$f_*\Q_{\sX}(-1)$ is constant, and concentrated in bidegree $(1,1)$.
Proposition \ref{6.4} implies that one has a decomposition
$$
R^{2}g_{2*}\Q_{\sZ_2}\oplus \W' \simeq
\bigoplus_{1}^{d-1} (R^1 g_{1*}\Q_{\sZ_1})_i \otimes
H^1(\Sigma_d,\Q_{\Sigma_d})_{d-i}\oplus \W,
$$
with $\W'$ and $\W$ constant and concentrated in degree $(1,1)$.
In particular, $R^{2}g_{2*}\Q_{\sZ_2}$ is a sub-variation of $\Q$-Hodge
structures of
$$R^1 g_{1*} \Q_{\sZ_1}\otimes H^1(\Sigma_d,\Q_{\Sigma_d})\oplus \W.
$$
Applying \ref{6.4} a second time, one finds
$$
R^{3}g_{3*}\Q_{\sZ_3}\simeq
\bigoplus_{1}^{d-1} (R^2 g_{2*} \Q_{\sZ_2})_i \otimes
H^1(\Sigma_d,\Q_{\Sigma_d})_{d-i}
\oplus \bigoplus^{d-1} R^{1}g_{1*}\Q_{\sZ_1}(-1),
$$
and $R^{3}g_{3*}\Q_{\sZ_3}$ is a sub-variation of $\Q$-Hodge
structures in
\begin{multline*}
R^1 g_{1*} \Q_{\sZ_1}\otimes H^1(\Sigma_d,\Q_{\Sigma_d})\otimes
H^1(\Sigma_d,\Q_{\Sigma_d})\oplus \W\otimes H^1(\Sigma_d,\Q_{\Sigma_d})\\
\oplus\bigoplus^{d-1}R^{1}g_{1*}\Q_{\sZ_1}(-1).
\end{multline*}
\begin{corollary}\label{7.5}
The set $S$ of points $s\in S$ for which $g_2^{-1}(s)$ and $g_3^{-1}(s)$
both have complex multiplication is dense in $S$.
\end{corollary}
\begin{proof}
Choose in \ref{7.4} a point $s\in S'$ such that $F_1=g_1^{-1}(s)$
has complex multiplication.
\begin{claim}\label{7.6}
$F_1\times\Sigma_5$,, $F_1\times \Sigma_5\times \Sigma_5$ and $\W\otimes H^1(\Sigma_5,\Q_{\Sigma_5})$
have all complex multiplication.
\end{claim}
\begin{proof}
It is well known that the Jacobian of every Fermat curve is of CM type
(apply for example the results in \cite{Shi}). Hence the first and the second part of \ref{7.6} follows
from \ref{7.1}, d), whereas the last one follows from \ref{7.1}, c) and d).
\end{proof}
Writing $F_3=g_3^{-1}(s)$ one has inclusions of polarized $\Q$-Hodge structures
$$
H^2(F_2,\Q_{F_2})\subset H^1(F_1,\Q_{F_1})\otimes H^1(\Sigma_5,\Q_{\Sigma_5})
$$
and
\begin{multline*}
H^3(F_3,\Q_{F_3})\subset H^1(F_1,\Q_{F_1})\otimes H^1(\Sigma_5,\Q_{\Sigma_5})\otimes
H^1(\Sigma_5,\Q_{\Sigma_5})\\
\oplus H^1(\Sigma_5,\Q_{\Sigma_5})\otimes \W\oplus H^1(F_1,\Q_{F_1})(-1).
\end{multline*}

By Lemma \ref{7.1}, a), b) and d), and by Claim \ref{7.6}
$$\Hg(H^1(F_1,\Q_{F_1})\otimes H^1(\Sigma_5,\Q_{\Sigma_5}))$$
and
\begin{multline*}
\Hg(H^1(F_1,\Q_{F_1})\otimes H^1(\Sigma_5,\Q_{\Sigma_5})\otimes
H^1(\Sigma_5,\Q_{\Sigma_5}) \oplus \\
H^1(\Sigma_5,\Q_{\Sigma_5})\otimes \W\oplus H^1(F_1,\Q_{F_1})(-1))
\end{multline*}
are both commutative.

By Deligne \cite{De1}, a polarized $\Q$-Hodge structure is semi simple.
Hence, the above inclusions of Hodge structures induce direct sum decompositions of
polarized $\Q$-Hodge structures. By the second part of Lemma \ref{7.1}, a), one has
surjective homomorphisms
$$
\Hg(H^1(F_1,\Q_{F_1})\otimes H^1(\Sigma_5,\Q_{\Sigma_5}))\>>> \Hg(H^2(F_2,\Q_{F_2}))
$$
and
\begin{multline*}
\Hg(H^1(F_1,\Q_{F_1})\otimes H^1(\Sigma_5,\Q_{\Sigma_5})\otimes
H^1(\Sigma_5,\Q_{\Sigma_5}))\times
\Hg(H^1(\Sigma_5,\Q_{\Sigma_5})\otimes \W)\times\\
\Hg(H^1(F_1,\Q_{F_1})(-1))\>>> \Hg(H^3(F_3,\Q_{F_3}))
\end{multline*}
Since the groups on the left hand sides are commutative, the groups on the right hand sides
are commutative, as well.
\end{proof}
\begin{proof}[Proof of Theorem \ref{0.2}]
Up to now, we only constructed a two-dimensional subscheme of $M_{5,4}$
with a dense set of CM-points. To get the second copy of $M_{5,1}$
we apply \ref{5.4}. There we constructed a generically finite morphism
$S\times S \to M_{5,4}$ which is induced by a family $g':\sZ' \to S\times S$.
Let us assume again, that $\Sigma$ is a curve with complex multiplication,
corresponding to a point $s \in S$. Restricting $g'$ to $S \times \{s\}$
one obtains a family $g_s:\sZ_s \to S$ satisfying the assumptions made in Section \ref{variation}
(for $\sZ_1$ instead of $\sX$). By Proposition \ref{6.4}
one finds an inclusion
$$ R^3g_{s*}\Q_{\sZ_s/S}\subset (R^2g_{2*}\Q_{\sZ_2/S})\otimes H^1(\Sigma,\Q_{\Sigma})
\oplus\bigoplus^4 R^1g_{1*}\Q_{\sZ_1/S}(-1).$$
The set of points in $s'\in S=S\times \{s\}$ where both, $H^2(F_2,\Q_{F_2})$
and $H^1(F_1,\Q_{F_1})$ have complex multiplication, is dense, and repeating
the argument used to prove \ref{7.5} one obtains that the $CM-$points are dense in the
image of $S\times S$ in $M_{5,4}.$    \\

The map $S\times S\to M_{5,4}$ is rigid. Otherwise one would get an extension
$$ g': \sZ'\>>> S\times S\times T,$$
such that $\dim T\geq 1$ and the induced map $S\times S\times T\to M_{5,4}$
is generically finite. Hence, by Proposition 3.1 $\dim S=\dim T=1$, a contradiction.
\end{proof}
\begin{remark}\label{7.7}
The construction in this Section is related to the one in \cite{FL}.
There it is shown that, starting from a universal family
$f:\sX\to S$ for $\sM_{5,1}$ and the corresponding second iterated $5$-fold covering
$g_3: \sZ_3\to S$, one obtains variations of Hodge structures
$\R^3g_{3*}\Q_{\sZ_3}$ which provide a uniformization of $M^{(3)}_{5,4}$
as a two-dimensional ball quotient.

Let $\V_1$ be the variation of $\Q$-Hodge structures of the $5$-fold covering
$g_1:\sZ_1\to S.$ By \cite{DM} $\V_1$ is the direct sum of
two unitary local systems and of two non-unitary local systems, dual to each
other. Moreover, by adding stable points Deligne and Mostow have defined  a
partial compactification $(M_{5,1})_{st}\supset M_{5,1}$  (see \cite{DM}, page 25 and Section 4),
such that each of the two non-unitary local systems induces a uniformization of $(M_{5,1})_{st}$
as a two dimensional ball quotient.

In Corollary \ref{7.4} the variation of $\Q$-Hodge structures $\V_3$, given by the family
$g_3:\sZ_3\to S$ is a direct factor of
$$\V_1\otimes H^1(\Sigma_5,\Q)^{\otimes 2}\oplus\W\otimes H^1(\Sigma_5,\Q)\oplus\V_1^{\oplus 4}.$$
This implies that every $\C$-irreducible non unitary direct factor of $\V_3$
provides a uniformization of $(M_{5,1})_{st}$ as a ball quotient.

Remark that the family $(g':\sZ'\to S\times S) \in \sM_{5,4}(S\times S)$ considered at the
end of the proof of \ref{0.2} gives the product of two $2$-dimensional ball quotients in a
partial compactification of $M_{5,4}$.
\end{remark}
%%%%%%%%%%%%%%%%%%%%%%%% References %%%%%%%%%%%%%%%%%
%\bibliographystyle{plain}

\end{document}